\documentclass[12pt]{article}
\usepackage{amsmath,amsthm,amssymb,color,graphicx}
\usepackage{subcaption}
\usepackage{url}

\begin{document}
\today\vspace{1.0cm}

\newtheorem{thm}{Theorem}[section]
\newtheorem{cor}[thm]{Corollary}
\newtheorem{lem}[thm]{Lemma}
\newtheorem{prop}[thm]{Proposition}
\newtheorem{rem}[thm]{Remark}
\newtheorem{defi}[thm]{Definition}
\newtheorem{conj}[thm]{Conjecture}

\renewcommand{\geq}{\geqslant} 
\renewcommand{\leq}{\leqslant} 
\renewcommand{\ge}{\geqslant} 
\renewcommand{\le}{\leqslant}

\newcommand{\gleb}[1]{\begin{center}\fbox{\parbox{.8\textwidth}{Gleb: #1}}\end{center}}
\newcommand{\tolia}[1]{\begin{center}\fbox{\parbox{.8\textwidth}{Anatoly: #1}}\end{center}}
\newcommand{\serg}[1]{\begin{center}\fbox{\parbox{.8\textwidth}{Sergey: #1}}\end{center}}
\newcommand{\maia}[1]{\begin{center}\fbox{\parbox{.8\textwidth}{Maia: #1}}\end{center}}

\begin{center}

{\bf\Large  Global Stability and Periodicity in a Glucose-Insulin Regulation Model with a Single Delay}

\vspace{8mm}

{\Large M. Angelova$^1$, G. Beliakov$^1$, A. Ivanov$^2$,  and S. Shelyag$^1$}




{\em $^1$ School of Information Technology, Deakin University Geelong, Australia}

{\em $^2$ Department of Mathematics, Pennsylvania State University, USA} \\









\end{center}

\vspace{8mm}
\begin{abstract}\noindent
A two-dimensional system of differential equations with delay modelling the glucose-insulin interaction processes in the human body is considered. Sufficient conditions are derived for the unique positive equilibrium in the system to be globally asymptotically stable. They are given in terms of the global attractivity of the fixed point in a limiting interval map. The existence of slowly oscillating periodic solutions is shown in the case when the equilibrium is unstable. The mathematical results are supported by extensive numerical simulations. 
It is deduced that typical behaviour in the system is the convergence to either a stable periodic solution or to the unique stable equilibrium. The coexistence of several periodic solutions together with the stable equilibrium is demonstrated as a possibility.
\end{abstract}


\section{Introduction}
\label{Intro}



This paper deals with further qualitative and numerical analyses of the Sturis-Bennett-Gourley 
model of the glucose-insulin interaction in the human body. The model was proposed in \cite{BenGou04-1} 
as a simple two-dimensional system of nonlinear differential equations with one delay. 
It can also be viewed as an abbreviated simplified version of a more complex differential model with 
two delays \cite{LiKua07,LiKuaMas06}.


Only a limited number of parameters of the physiologically closed glucose-insulin interaction system are accessible for direct measurements. Therefore, mathematical modelling is required to facilitate the estimation of the narrow physiological range of the glucose-insulin system components \cite{marchetti2016}. The physiological time delay plays an important part in the regulation and the feedback of the system and determines the important properties of the mathematical description of the model.

The glucose-insulin regulatory system is a key component in the metabolism of the human body. The pancreas and the liver regulate the production of insulin and glucose respectively in order to maintain normal level of the blood glucose. Failure to do this can result in high blood sugar and  diabetes which is related to many other long term health problems. 
Within this regulation, both rapid (period  $\sim 6-15$ mins) and ultradian (period $\sim 80-180$ mins) oscillations of insulin  have been observed \cite{Satin2015}, along with glucose oscillations (period $\sim 80-150$ mins) \cite{Scheen96}. The ultradian oscillations were  first discovered by Hansen \cite{Hansen23} and observed during fasting, meal ingestion, continuous enteral nutrition and constant glucose infusion. 

Due to the monotone nature of the nonlinearities involved in the two-dimensional system it is a natural 
conjecture that some degree of  the simplicity in the dynamical behaviour of the solutions should be observed 
in the model. 
Such simpler dynamics would be consistent with known behaviours in scalar differential delay equations 
with monotone nonlinear feedback \cite{M-PSel96b,M-PWal94,Wal81a}. 
In fact, some of the known theoretical results can be extended to the Sturis-Bennett-Gourley model (see further details in the next section Preliminaries).

The primary objective of this paper is to derive conditions for the global asymptotic stability 
in the model as well as conditions for the existence of nontrivial periodic  solutions. 
We apply and further develop some of our earlier approaches and results to study differential delay equations and systems via underlying finite-dimensional discrete maps \cite{BraHasIvaTro20,HalIva93,IvaBLW19,IvaSha91}. In our analyses in this paper we also use prior related results, in  particular those obtained in \cite{BenGou04-1,BenGou04-2,BenGou04-3,M-PSel96a,M-PSel96b,M-PWal94}. 
The differential delay model we consider was introduced and studied in paper \cite{BenGou04-1}. 
We further investigate its properties and interpret some of them from an alternative prospective via a limiting interval map. The principal result of paper \cite{M-PSel96b}, an extension of the Poincar\'{e}-Bendixson Theorem
\cite{Poincare,Bendixson, Coddington}
to delay differential equations, can be applied to the Sturis-Bennett-Gourley system as well. As a result one can conclude that its every solution converges to either the unique equilibrium or to a periodic solution. In this paper, we demonstrate various possibilities of such eventual dynamics through several numerical examples.



The approach we propose to analyse the system is through a limiting one-dimensional map, which is formally obtained from the original differential delay system when the delay goes to infinity. The resulting interval map is relatively simple: it is given by a monotone decreasing function and can only have a unique fixed point and cycles of period two  (which can be either attracting, or repelling, or a combination of the two possibilities).

The dynamics of the limiting interval map largely determine the dynamics of our two-dimensional differential delay system. When the only fixed point of the interval map is globally attracting then the corresponding unique steady state of the delay system is also globally asymptotically stable (for arbitrary delay $\tau>0$). However, when the only fixed point of the interval map is repelling the dynamics in the differential delay system are varied and dependent on the size of delay $\tau$. When the delay is small enough, $0\le\tau<\tau_0$ for some $\tau_0>0$, the unique equilibrium of the differential delay system is locally asymptotically stable. Note that only the local stability of the equilibrium can be claimed, as examples are possible that the delay system has stable periodic solutions away from the stable constant equilibrium (see an example in Section 4). 
When the delay becomes large enough, $\tau>\tau_1$ for some $\tau_1>0$, then the differential delay system has a slowly oscillating periodic solution. In this case the corresponding characteristic equation of the linearised system about the steady state has a leading pair of complex conjugate solutions $\alpha_0\pm\beta_0 i$ with the positive real part $\alpha_0>0$ and the imaginary part within the interval $\beta_0\in(0,\pi/\tau)$.  This leading eigenvalue makes the slow oscillation in the system typical, in agreement with known results for the case of similar scalar differential delay equations \cite{M-PWal94,Wal81a}. 

It is important to note that from the  mathematical point of view there is no uniqueness in the model for stable periodic solutions or for the stable equilibrium. We construct explicit examples of our two-dimensional differential delay system when the coexistence of two stable periodic solutions is observed (Section 4). We also demonstrate the possibility when a stable periodic solution coexists together with the locally stable equilibrium (Section 4). The examples are easily generalised to the case when any finite number of stable periodic solutions can coexist with or without the locally attracting equilibrium. This multi-stability phenomenon implies the utmost importance of the proper choices of the nonlinearities when the differential delay system is suggested as an actual model of a particular applied problem. For some of such known and available in the literature models we numerically observe the uniqueness of the stable periodic solution and its global attractivity within the admissible set of initial conditions. 

The novelty of our approach is that we derive a simple one-dimensional dynamical system (interval map) and use it to determine the global dynamical properties of our infinite-dimensional dynamical system with delay. The typical dynamical behaviours in our model are also simple - it is the convergence to either a stable periodic solution or to the unique stable equilibrium, or a combination of such behaviours.

The paper is organised as follows. Section 2 describes the foundations of the mathematical problem. 
It includes references to and a brief review of existing closely related results obtained by others.
In Section 3, we derive the main analytical results, which are then numerically confirmed in Section 4.  
A concluding summary and a brief discussion are given in Section 5.

\section{Preliminaries}\label{Prelim}

\subsection{Differential Delay Model and Assumptions}\label{DDMA}

Consider the system of differential equations with delay \cite{BenGou04-1}, 
\begin{eqnarray}\label{MS}
            I^\prime(t) &=& f_1(G(t))-\frac{1}{\tau_0} I(t) \\
  \nonumber G^\prime(t) &=& G_{in}-f_2(G(t)))-qG(t)f_4(I(t))+f_5(I(t-\tau)),
\end{eqnarray}
where $I$ and $G$ represent the relative concentrations  of insulin and glucose, respectively, ${1}/{\tau_0}$ is the insulin degradation rate, $G_{in}$ is the external glucose input. The function 
$f_1$ corresponds to pancreatic insulin production, dependent on glucose concentration, and $f_2$ is the glucose consumption by the brain. The third term in the second equation represents the insulin-dependent glucose utilisation in the muscles, while the last term, $f_5$, represents the hepatic glucose production.  $\tau$ is the time delay between plasma insulin production and its effect on hepatic glucose production. 

The system is considered under the following assumptions:
\begin{itemize}
  \item[(H1)] Functions $f_1(u), f_2(u), f_4(u), f_5(u)$ are non-negative and continuously differentiable for $u\ge0$, with $f_3$ defined as $f_3(u)=q u$ for convenience.
  Real parameters $\tau_0, G_{in}, q, \tau$ are all positive;
  \item[(H2)] $f_1(u)>0, f_1^\prime(u)>0, \forall u>0, f_1(0)=a_0>0$ and $\lim_{u\to\infty}f_1(u)=a>0$;
  \item[(H3)] $f_2(u)>0, f_2^\prime(u)>0, \forall u>0, f_2(0)=0$ and $\lim_{u\to\infty}f_2(u)=b>0$;
  \item[(H4)] $f_4(u)>0, f_4^\prime(u)>0, \forall u>0, f_4(0)=d>0$ and $\lim_{u\to\infty}f_4(u)=e>0$;
  \item[(H5)] $f_5(u)>0, f_5^\prime(u)<0, \forall u>0, f_5(0)=h>0$ and $\lim_{u\to\infty}f_5(u)=0$.
\end{itemize}
The assumptions (H1)-(H5) are derived from and justified by the physiological mechanisms of the glucose-insulin interaction in the human body, see e.g. papers \cite{BenGou04-1,BenGou04-2,LiKua07,LiKuaMas06} for additional details.

The phase space of system (\ref{MS}) is defined as $\mathbb X= C([-\tau,0], \mathbb R_+)\times \mathbb R_+$ where $\mathbb R_+:=\{x\in\mathbb R\vert x\ge0\}$.
For arbitrary initial function $\psi=(\varphi(s), u)\in\mathbb X$ the corresponding solution $\mathbf x=\mathbf x(t,\psi)=(I(t),G(t))$ to system (\ref{MS}) can be constructed by the standard step methods \cite{BelCoo63,DieSvGSVLWal95,HalSVL93}.
We assume that such solutions exist for an arbitrary initial function $\psi\in\mathbb X$ and all $t\ge0$ (which is the case under the assumptions that nonlinearities $f_1, f_2, f_4$ are continuously differentiable).

It is an easy observation that positive initial data for system (\ref{MS}) result in solutions that are positive 
for all $t\ge0$. More precisely, if the initial function $\psi=(\varphi(s),u), s\in[-\tau,0],$ is such that 
$u>0, \varphi(s)\ge0\; \forall s\in[-\tau,0]$ and $\varphi(0)>0$ then $I(t)>0$ and $G(t)>0$ holds for all 
$t\ge0$ (see further details and proof in \cite{BenGou04-1}).
It can also be shown that both components $I$ and $G$ of all solutions to system (\ref{MS}) are bounded from above and bounded away from zero. Moreover, a stronger property called the persistence can be established here. It says that positive constants $m_I, m_G$ and $M_I, M_G$ can be identified, independent of particular initial data, such that for an arbitrary initial function $\psi\in\mathbb X$ and the corresponding solution 
$\mathbf x(t,\psi)=(I(t), G(t))$ to system (\ref{MS}) there exists a time moment $T=T(\psi)$ such that the following holds
\begin{equation}\label{Persist}
  0<m_I\le I(t)\le M_I<\infty,\;  0<m_G\le G(t)\le M_G<\infty, \quad \text{for all}\quad t\ge T.
\end{equation}

These and other basic properties of the solutions are proved in \cite{BenGou04-1} as Propositions 2.1, 2.2, and 2.4.
We will revisit them later in the paper from a different point of view.

\subsection{Linearization and Characteristic Equation}\label{LCE}
In this subsection we present well-known facts about the unique positive equilibrium of system (\ref{MS}), the linearized system about the equilibrium, and the characteristic equation of the linear system. More related details can be found in papers \cite{BenGou04-1,BenGou04-2} and \cite{LiKua07}.

Equilibria of differential delay system (\ref{MS}) are found by solving the nonlinear system
\begin{equation}\label{equil}
    I=\tau_0 f_1(G),\quad f_2(G)+qGf_4(I)=G_{in}+f_5(I),
\end{equation}
which reduces to a single scalar equation for $G$: $f_2(G)+qGf_4(\tau_0 f_1(G))=G_{in}+f_5(\tau_0 f_1(G))$.
It is straightforward to see that the latter has a unique positive solution $G_*>0$, implying that the original system  (\ref{MS}) has a unique equilibrium
$(I_*, G_*)$ where $I_*=\tau_0 f_1(G_*)>0$.

The linearized system about the positive equilibrium $(I_*,G_*)$ has the form 
\begin{eqnarray}\label{LMS}
            u^\prime(t) &=& -\frac{1}{\tau_0} u(t)+f_1^\prime(G_*) v(t) \\
  \nonumber v^\prime(t) &=& -[f_2^\prime (G_*)+q f_4(I_*)] v(t)-qG_*f_4^\prime (I_*) u(t)+f_5^\prime (I_*) u(t-\tau).
\end{eqnarray}
Note that system (\ref{LMS}) is also the linearization of the translated system (\ref{MS-0}) (derived in subsection \ref{TZE}). 
The characteristic equation of the linear system (\ref{LMS}) has the form
\begin{equation}\label{ChEq}
    (\lambda+\mu_1)(\lambda+\mu_2)+b+a\exp\{-\tau\lambda\}=0,
\end{equation}
where $\mu_1=1/\tau_0, \mu_2=f_2^\prime(G_*)+q f_4(I_*), b=q G_* f_1^\prime(G_*) f_4^\prime(I_*), a=-f_1^\prime(G_*)f_5^\prime(I_*).$ 
Since $f_1^\prime(G_*)>0, f_2^\prime(G_*)>0, f_4^\prime(I_*)>0$ and $f_5^\prime(I_*)<0$ then $\mu_1>0, \mu_2>0, b>0, a>0.$

The form (\ref{ChEq}) of the characteristic equation allows us to use known facts about its properties derived elsewhere, 
see e.g. \cite{adH79a,BraHasIvaTro20,LiKua07}.
The stability/instability of the zero solution of system (\ref{LMS}) is determined by the location of the solutions of the 
characteristic equation (\ref{ChEq}) in the complex plane. 
If all solutions of the characteristic equation have negative real parts (or are negative themselves) then the zero solution of (\ref{LMS}) 
is asymptotically stable.
If the characteristic equation (\ref{ChEq}) has a complex conjugate solution $\lambda=\alpha+i \beta$ with the positive real part $\alpha>0$ 
then the zero solution of (\ref{LMS}) is unstable.
In the latter case there exists the so-called leading pair of complex conjugate solutions $\lambda=\alpha_0\pm i \beta_0$ of the characteristic 
equation (\ref{ChEq}),  where $0<\beta_0<\pi/\tau$ and $\alpha_0>0$. The leading means that $\alpha_0>0$ is the largest real part among all solutions 
of (\ref{ChEq}). 
All other complex conjugate solutions $\lambda=\alpha_k\pm i \beta_k, k\in\mathbb N,$ of (\ref{ChEq}) satisfy $\alpha_0>\alpha_1>\alpha_2> \dots $ and $\beta_k\in[2k\pi/\tau, (2k+1)\pi/\tau]$. 
See Lemma 1 of \cite{adH79a} and Lemma 3 of \cite{BraHasIvaTro20} for more details and proofs.

The described above stability or instability of the zero solution of the linear system (\ref{LMS}) in terms of the eigenvalues of the characteristic equation (\ref{ChEq}) carry over to the nonlinear system (\ref{MS}). The constant solution $(I_*,G_*)$ of the latter has the same type of stability as $(0,0)$ of the linear system provided all the nonlinearities $f_1, f_2, f_4, f_5$ are $C^1$-smooth in a neighborhood of $(I_*,G_*)$ (see \cite{HirSma74}, Chapter 9, for further details and proof).

One of the main questions of interest in this paper is about the typical behaviour of solutions in system (\ref{MS}). In view of the strong monotonicity properties of the nonlinearities $f_1, f_2, f_4$ and $f_5$ some of the previously obtained results are applicable and can be used in our analysis.  The main result of paper \cite{M-PSel96b}, Theorem 2.1, implies that for every initial function $\psi\in\mathbb X$ the corresponding solution $\mathbf x(t,\psi)$ converges as $t\to\infty$ either to a periodic solution or to the unique equilibrium $(I_*,G_*)$ of system (\ref{MS}). The follow-up question is which of these possible behaviours are typical, and which ones can be observed in numerical simulations. It turns out that both stable periodic solutions or a stable equilibrium or a simultaneous coexistence of both are typical (as shown numerically in section \ref{numerical}).  It is also evident that there exist unstable periodic solutions (e.g. those separating neighboring stable periodic solutions; they usually cannot be observed numerically).

Oscillatory solutions in differential delay systems are typical in general. It is a known fact that all solutions to system (\ref{MS-0}) oscillate when the characteristic equation (\ref{ChEq}) has no real eigenvalues \cite{BraHasIvaTro20} (also see Proposition \ref{osc} of subsection \ref{PSs}). Due to the overall negative feedback  the slow oscillation can be typical in system (\ref{MS}) (see the related definition 3.9 in subsection \ref{osc}). In particular, initial functions $\psi=(\varphi,u)\in\mathbb X$ with $\varphi(s)-I_*>0\; \forall s\in[-\tau,0], u>0$ give rise to slowly oscillating solutions. This oscillation is typical as small perturbations of such initial functions leave them within the same initial set.

Another possibility for a typical behaviour is that solutions converge monotonically to the equilibrium $(I_*,G_*)$ as $t\to\infty$. This is the case when the characteristic equation (\ref{ChEq}) has real eigenvalues (which are then necessarily negative). Such a case implies the existence of the exponential decaying solutions to the linear equation (\ref{LMS}), and also the existence of solutions close to the exponential ones for the nonlinear system (\ref{MS}).

The results of papers \cite{M-P88,M-PWal94,Wal81a} for scalar equations suggest that typical behaviour of the oscillating solutions in system (\ref{MS}) is the eventual slow oscillation. That is, for almost all initial functions $\psi\in\mathbb X$ the corresponding solution $\mathbf x(t,\psi)$ is slowly oscillating for $t\ge T$ for some $T(\psi)\ge0.$ This means that if an initial function $\psi_0=(\varphi_0,u_0)\in\mathbb X$ is such that the solution $\mathbf x(t,\psi_0)$ is not slowly oscillating for $t\ge0$ then its every neighborhood $U(\psi_0)$ contains an initial function $\psi_1\in U_1$ such that the corresponding solution $\mathbf x(t,\psi_1)$ is eventually slowly oscillating.

Note that the problem of typical behaviour in a differential delay system is a very challenging one in general (many aspects of this problem remain unsolved even for the scalar case of a simple single equation; the basic paper \cite{M-PWal94} on the issue exists in a preprint form only). Therefore, we are in a position in this paper to only numerically verify the assumed theoretical results about the typical behaviours. The rigorous mathematical proofs will hopefully be accomplished  at a later time.

\subsection{Related Interval Maps}
\label{IMs}

In this subsection we recall some basic notions and definitions on interval maps related to the needs of this paper. Comprehensive expositions on the theory of one-dimensional maps can be found e.g. in monographs \cite{deMStr93,ShaKolSivFed97}.

Given a continuous map $F: L\rightarrow L$ of a closed interval $L\subseteq\mathbb R$ into itself a 
forward trajectory through an initial point $x_0\in L$ is defined as the set $\{F^n(x_0), n\in\mathbb N_0\}$
where $F^n=F\circ F\circ\dots\circ F$ is the $n^{\text{th}}$ iteration of 
$F$ ($F^0(x):=x;\, \mathbb N_0:=\mathbb N \cup\{0\}$).
A set $J\subset L$ will be called {\it invariant} under $F$ if $F(J)\subseteq J$. 
Note that a proper inclusion is allowed under this definition.

\begin{defi} (i) A fixed point $x=x_*$ of a continuous map $F$ of an interval $L\subseteq\mathbb R$ into itself is called ${\it attracting}$ if there exists an open interval $J\subseteq L$ such that $x_*\in J$, $f(J)\subseteq J$, and for every point $x\in J$ one has that $\lim_{n\to\infty} F^n(x)=x_*$ holds.\newline
\noindent
(ii) The largest connected interval $J\subseteq L$ with this property is called the  {\rm domain of immediate attraction}  of the fixed point $x_*$.
\noindent
(iii) A point $x_0$ is called periodic with period $m$ if $F^m(x_0)=x_0$ and $F^k(x_0)\ne x_0$ for every $1\le k\le m-1$. 
The corresponding set $\{x_0, x_1, \dots, x_{m-1}\}:=C_m$ is called a cycle of period $m$.
\end{defi}

Clearly that every point of the cycle $x_k\in C_m$ is periodic of period $m$ for the map $F$; it is also a fixed point for the map $F^m$.

The following statement is a well-known simple fact in the theory of interval maps. Its proof easily follows from related facts of Section 2.4 in \cite{ShaKolSivFed97}.

\begin{prop}\label{prop1}
For an arbitrary point $x\in J$ in the domain of immediate attraction of the fixed point $x_*$ there always exists a closed finite interval $L_0=L_0(x)\subset J$ such that $x\in L_0, F(L_0)\subseteq L_0$, and $\cap_{n\ge0}F^n(L_0)=x_*.$
\end{prop}

\begin{defi}
Let $x_*$ be an attracting fixed point of a continuous map $F$. An infinite set of intervals $\{L_n, n\in\mathbb N_0\}$ will be called a {\em squeezing sequence of imbedded intervals} if the following holds:
$$
L_{k+1}\subseteq L_{k}, F(L_k)\subseteq L_{k+1},\; {\rm and}\; \cap_{k\ge0}L_k=x_*.
$$
\end{defi}

It is evident that the sequence of intervals $L_k=F^k(L_0), n\in\mathbb N_0,$ in Proposition \ref{prop1} is a squeezing  imbedded sequence.
Given an initial point $x_0$ in the domain of immediate attracting it is also clear that a squeezing imbedded sequence of intervals containing its iterations always exists but is not uniquely defined in general.

\subsection{Translation to Zero Equilibrium}\label{TZE}

It was demonstrated in \cite{BenGou04-1,BenGou04-2} that under the assumptions (H1)-(H5) the system (\ref{MS}) 
has  unique equilibrium $(I_*,G_*), I_*>0, G_*>0,$ where $I_*=\tau_0 f_1(G_*)$ and $G_*$ is a unique positive
solution of the nonlinear equation $f_2(G)+qG f_4 (\tau_0 f_1(G))=G_{in}+f_5(\tau_0 f_1(G))$. 
As a matter of convenience, for various theoretical considerations and computational tasks of this paper 
it is advantageous to have this equilibrium shifted to the zero equilibrium state $(I_*,G_*)=(0,0)$.
One of the reasons for this need is that the equilibrium $(I_*,G_*)$ depends on all the parameters and functions
involved in system (\ref{MS}). 
Such shift is achieved by the change of the dependent variables by
\begin{equation}\label{ChVar}
    x(t)=I(t)-I_*,\quad y(t)=G(t)-G_*.
\end{equation}
System (\ref{MS}) is then transformed into the following one:
\begin{eqnarray}\label{MS-0}
            x^\prime(t) &=& F_1(y(t))-\frac{1}{\tau_0} x(t) \\
  \nonumber y^\prime(t) &=& -F_2(y(t))-q\,f_4(I_*+x(t))\, y(t)-q\,G_*\, F_4(x(t))+F_5(x(t-\tau)),
\end{eqnarray}
where 
$F_1(y)=f_1(y+G_*)-f_1(G_*), F_2(y)=f_2(y+G_*)-f_2(G_*), F_4(x)=f_4(x+I_*)-f_4(I_*), F_5(x)=f_5(x+I_*)-f_5(I_*).$
Functions $F_1, F_2, F_4$ are strictly monotone increasing and satisfying the positive feedback condition $y\cdot F_i(y)>0$ for $y\ne0, i=1,2,4$.
Function $F_5(x)$ is strictly decreasing and satisfying the negative feedback assumption $x\cdot F_5(x)<0$ for $x\ne0$.
System (\ref{MS-0}) has the unique zero equilibrium $(x,y)=(0,0)$.

We will perform most of our numerical simulations for systems of type (\ref{MS-0}). By using appropriate inverse transformations such systems can always be represented in the form of the original system (\ref{MS}).


\section{Main Results}
\label{main_results}

\subsection{Limiting Interval Map}
\label{lim_int_map}
In this sub-section we derive a limiting interval map for the differential delay system (\ref{MS}) as $\tau\to\infty$. First we transform system (\ref{MS}) to one with the normalised delay $\tau=1$ by rescaling the independent variable by $t=\tau\cdot s$. It is a straightforward calculation then that reduces system (\ref{MS}) to the following form
\begin{eqnarray}\label{MS1}
            \frac{1}{\tau}\,I^\prime(s) &=& f_1(G(s))-\frac{1}{\tau_0} I(s) \\
  \nonumber \frac{1}{\tau}\,G^\prime(s) &=& G_{in}-f_2(G(s))-qG(s)f_4(I(s))+f_5(I(s-1)).
\end{eqnarray}
By taking the limit as $\tau\to\infty$ the latter becomes a system of functional difference equations:
\begin{equation}\label{SDE}
I(s)=\tau_0\,f_1(G(s)),\quad f_2(G(s)) + qG(s)f_4(I(s)) = G_{in}+f_5(I(s-1)),
\end{equation}
which in turn is further reduced to a single scalar difference equation for the variable $G$:
\begin{equation}\label{DE}
f_2(G(s)) + qG(s)f_4(\tau_0f_1(G(s))) = G_{in}+f_5(\tau_0f_1(G(s-1))).
\end{equation}
It is easy to see, based on the assumptions (H1)-(H5), that the function $F$ in the left hand side of equation (\ref{DE}), $F(G):=f_2(G)+qG\,f_4(\tau_0f_1(G))$, satisfies:
\begin{equation}\label{F}
F(0)=0, F^\prime(G)>0,~\forall G\ge0,\quad \text{and}\quad \lim_{G\to\infty} F(G)=\infty.
\end{equation}
Likewise, the function $H$ in the right hand side of equation (\ref{DE}), $H(G):=G_{in}+f_5(\tau_0f_1(G))$, satisfies:
\begin{equation}\label{H}
H(0)=H_0>0, H^\prime(G)<0,~\forall G\ge0,\quad \text{and}\quad \lim_{G\to\infty} H(G)=H_{\infty}>0.
\end{equation}
Therefore, the inverse function $F^{-1}$ exists, and equation (\ref{DE}) can be explicitly solved for $G(s)$ as follows:
\begin{equation}\label{IM}
G(s)=F^{-1}(H(G(s-1)))=:\Phi(G(s-1)),
\end{equation}
where the composite function $\Phi=F^{-1}\circ H$ is defined and continuous on $\mathbb{R_+}=\{G\vert G\ge0\}$.
Besides, due to assumptions $(H1)-(H5)$, function $\Phi(\cdot)$ is continuously differentiable on $\mathbb R_{+}$ with
\begin{equation}\label{Phi}
\Phi^{\prime}(u)<0 \:\; \forall u\in\mathbb R_+\quad \text{and}\quad 
\lim_{u\to0}\Phi(u)=\Phi_0>0, \lim_{u\to\infty}\Phi(u)=\Phi_{\infty}>0.
\end{equation}
The values $\Phi_0, \Phi_{\infty}$ are easily calculated as:
\begin{equation}\label{Phi-lim}
\Phi_0=F^{-1}(G_{in}+f_5(\tau_0a_0)),\qquad \Phi_{\infty}=F^{-1}(G_{in}+f_5(\tau_0a)).
\end{equation}
The asymptotic properties of solutions of equation (\ref{IM}) are completely determined by the dynamical properties of the iterations of the interval map $\Phi$.
A comprehensive theory of such equations is given in the monograph \cite{ShaMaiRom93}. All relevant properties on interval maps can be found in monographs \cite{deMStr93,ShaKolSivFed97}.

A convenient look at system (\ref{SDE}) and equation (\ref{IM}) is via difference equation notations. By denoting $G(t):=G_n, I(t):=I_n, I(t-1):=I_{n-1}, n\in\mathbb N$ system (\ref{SDE}) is rewritten as
$$
f_2(G_n)+q\,G_n\,f_4(I_n) = G_{in}+f_5(I_{n-1}),\quad I_n=\tau_0 f_1(G_n).
$$
The difference equation (\ref{IM}) is represented then as $G_n=\Phi(G_{n-1}), n\in\mathbb{N}$.


\subsection{Principal Results}
\label{PRs}

Based on property (\ref{Phi}) we can build a sequence of imbedded intervals for map $\Phi$ as follows. 
Set $L_0:=\mathbb R_+$ and $L_1:=\Phi(L_0)=\Phi(\mathbb R)=[\Phi_\infty,\Phi_0]\subset L_0.$ 
Proceed then recursively as $L_2=\Phi(L_1)\subset L_1, \dots, L_{n+1}=\Phi(L_n)\subset L_n,n\in\mathbb N_0.$
Define the limiting set $L_*$ by $L_*:=\cap_{n\ge0}\, L_n:=[\alpha_*,\beta_*]$.

The set $L_*$ is either a single point or a closed interval with a non-empty interior. In the first case one has that $\alpha_*=\beta_*=G_*$. 
In the second case the  endpoints $\{\alpha_*,\beta_*\}$ form a cycle of period two: $\alpha_*=\Phi(\beta_*), \beta_*=\Phi(\alpha_*), \Phi(L_*)=L_*, G_*\in {int}\,(L_*)$.

The sequence $\{L_n\}$ of imbedded interval for the component $G$ generates the sequence of imbedded intervals $\{J_n\}$ for the component $I$ 
through the first difference equation of system (\ref{SDE}) by $J_n:=\tau_0f_1(L_n), n\in\mathbb N_0$. 

We shall formally distinguish the following subcases for the set $L_*$ and its structure:
\begin{itemize}
    \item[$(A_1)$] 
    The set $L_*$ is a single point $G_*$. It is the only fixed point of the interval map $\Phi$ which is then globally attracting on $\mathbb R_+$: 
    for every initial point $G_0\in\mathbb R_+$ one has $\lim_{n\to\infty}\Phi^n(G_0)=G_*$;
    \item[$(A_2)$]
    The set $L_*$ is a closed non-empty interval $[\alpha_*,\beta_*], \alpha_*\ne\beta_*,$.
    Then $\Phi(L_*)=L_*,$ and $\{\alpha_*,\beta_*\}$ is a cycle of period two, $\alpha_*=\Phi(\beta_*), \beta_*=\Phi(\alpha_*),$ with $G_*\in {int}\,(L_*)$. 
    Assume also that the two-cycle is globally attracting: for every initial point $G_0\in\mathbb R_+, G_0\ne G_*,$ its forward iterations converge to the cycle: 
    $\Phi^n(G_0)\longrightarrow \{\alpha_*,\beta_*\}$ as $n\to\infty$.
    In addition, it is assumed that the generic condition $\Phi^\prime(G_*)<-1$ holds;
    \item[$(A_3)$] 
    The set $L_*$ is an interval formed by a cycle $\{\alpha_*,\beta_*\}$ of period two which is locally attracting only. In addition, it is assumed that the condition $\Phi^\prime(G_*)<-1$ holds;
    \item[$(A_4)$]
    The set $L_*$ is an interval formed by a cycle $\{\alpha_*,\beta_*\}$ of period two. In addition, the fixed point $G_*$ is locally attracting: 
    there exists an interval $(\gamma_*,\delta_*)$ such that  for every initial point $G_0\in(\gamma_*,\delta_*)$ one has $\lim_{n\to\infty}\Phi^n(G_0)=G_*$;
\end{itemize}

Below we state the principal results of this paper which are essentially dependent and built based on the structure of the limiting sets 
of the map $\Phi$ as described above by properties $(A_1)-(A_4)$.

\begin{thm}\label{ThmA1}
Suppose that in addition to $(H_1)-(H_5)$ assumption $(A_1)$ holds. Then for every delay $\tau>0$ the unique equilibrium $(I_*,G_*)$ of system (\ref{MS}) is globally asymptotically stable: 
for an arbitrary initial data $\psi=(\varphi,u)\in\mathbb X$ the corresponding solution $(I(t),G(t))$ satisfies: $\lim_{t\to\infty}G(t)=G_*, \lim_{t\to\infty}I(t)=I_*$.
\end{thm}

Theorem \ref{ThmA1} is a strong delay independent result about the global asymptotic stability of the unique equilibrium $(I_*,G_*)$ of system
(\ref{MS}) with infinite-dimensional phase space based on the global attractivity of the corresponding fixed point in a simple limiting 
interval map defined by the real-valued function $\Phi$. Theorem \ref{ThmA1} is proved in subsection \ref{IPGAS} (as Theorem \ref{GAS}). 

\medskip\noindent
Note that a closely related result to Theorem \ref{ThmA1} is proved in paper \cite{BenGou04-1} 
as Theorem 3.2. 
However, their approach and method of proof are different from what we use in present paper. 
We reflect more on paper \cite{BenGou04-1} in the Conclusion section.

\begin{thm}\label{ThmA2}
Suppose that in addition to $(H_1)-(H_5)$  assumption $(A_2)$ holds. Then there exists $\tau_0>0$ such that for every delay $\tau>\tau_0$ system (\ref{MS}) has a slowly oscillating periodic solution.
\end{thm}

The slow oscillation of solutions here means that both components $I(t)$ and $G(t)$ are slowly oscillating functions about their respective equilibrium values $I_*$ and $G_*$. 
Therefore, $I(t)-I_*$ and $G(t)-G_*$ are slowly oscillating functions with their successive zeros separated by a time span larger than the delay $\tau$.
For more complete definitions and statements see details in Subsection \ref{PSs}.

Under assumption $(A_2)$ the cycle $\{\alpha_*,\beta_*\}$ is globally attracting on $\mathbb R_+$: for every initial value $G_0\in\mathbb R_+, G_0\ne G_*,$ the sequence of its iterations $\Phi^n(G_0)$ is attracted by the cycle as $n\to\infty$. 
This means that both sequences $\Phi^{2n}(G_0)$ and $\Phi^{2n+1}(G_0), n\in\mathbb N_0,$ are monotone and converge to either $\alpha_*$ or $\beta_*$ (depending on the location of $G_0$ in relation to the fixed point $G_*$).
To show the existence of a slowly oscillating periodic solution to system (\ref{MS}) we use the standard and well developed techniques of the ejective fixed point theory \cite{DieSvGSVLWal95,HalSVL93}.   
To that end the main points we have to show holding true for system (\ref{MS}) are:
\begin{itemize}
    \item[(i)] 
    Construction of a cone of initial data for system (\ref{MS}) and a non-linear map which maps the cone into itself. The map usually is an appropriately defined shift operator along the solutions;
    \item[(ii)] 
    Existence of a leading eigenvalue to the characteristic equation (\ref{ChEq}) with the largest positive real part and the imaginary part within the range $(0,\pi/\tau)$;
    \item[(iii)] 
    The compactness of the shift operator along solutions of the system starting on the cone.
\end{itemize}
The outline of the proof of the existence of periodic solutions is given in Subsection \ref{PSs}.

Note that in general the ejective fixed point techniques do not address the issue of the uniqueness of the slowly oscillating periodic solution. 
The existence of periodic solutions can only be proved; the periodic solutions can be non-unique in many cases.
This is true for all the classes of delay equations and systems to which they were applied, including our system (\ref{MS}). 
However, the uniqueness of the globally attracting cycle $\{\alpha_*,\beta_*\}$ of period two for the one dimensional map $\Phi$ seems to yield the uniqueness of a stable slowly oscillating periodic solution to system (\ref{MS}).
This fact can be verified numerically. We have done it for two classes of the nonlinearities $f_i, i=1,2,4,5,$ used in applications: Hill type functions \cite{HuaEasAng15} and exponential functions \cite{LiKuaMas06}. 

\begin{thm}\label{ThmA3}
Suppose that in addition to $(H_1)-(H_5)$  assumption $(A_3)$ holds. 
Then there exist multiple choices of the nonlinearities $f_1,f_2,f_4,f_5$ and of the parameter values $\tau_0,\tau, q$ such that 
system (\ref{MS}) possesses at least two slowly oscillating periodic solutions.
\end{thm}

The key assumption in $(A_3)$ is that the two-cycle is locally attracting only; therefore, there exists another cycle of period two.
Since $\Phi^\prime(G_*)<-1$ the fixed point $G_*$ is repelling. 
Hence, there exists the minimal cycle of period two, $\{\gamma_*,\delta_*\}$, such that the open interval $(\gamma_*,\delta_*)\ni G_*$ is attracted to it.
In addition, the inequalities $\alpha_*<\gamma_*<G_*<\delta_*<\beta_*$ hold. Both cycles $\{\alpha_*,\beta_*\}$ and $\{\gamma_*,\delta_*\}$ are at least one-sided attracting.
The structure of the map $\Phi$ on the set $[\alpha_*,\gamma_*]\cup[\delta_*,\beta_*]$ can be arbitrary; however, since $\Phi$ is monotone decreasing, it can only contains additional cycles of period two.

The non-uniqueness of the two-cycle $\{\alpha_*,\beta_*\}$, and the existence of the second two-cycle $\{\gamma_*,\delta_*\}$, seem to be an important factor for the existence of multiple periodic solutions to system (\ref{MS}). 
We construct an example of system (\ref{MS}) when it has one slowly oscillating periodic solution related to the smallest two-cycle  $\{\gamma_*,\delta_*\}$ and the second slowly oscillating periodic solution related to the largest two-cycle $\{\alpha_*,\beta_*\}$ (see an example in subsection \ref{MPS}). Such example can be easily generalized to produce any finite number of slowly oscillating periodic solutions to system (\ref{MS}).

\begin{thm}\label{ThmA4}
Suppose that in addition to $(H_1)-(H_5)$ assumption $(A_4)$ holds. 
Then there exist multiple choices of the nonlinearities $f_1,f_2,f_4,f_5$ and of the parameter values $\tau_0,\tau, q$ such that 
system (\ref{MS}) possesses both a slowly oscillating periodic solution and the locally attracting equilibrium $(I_*,G_*)$.
\end{thm}

The principal difference between Theorem \ref{ThmA4} and Theorem \ref{ThmA3} is that the fixed point $G_*$ is attracting for the map $\Phi$ in the latter (while it was repelling for the former).
Therefore, its minimal two-cycle  $\{\gamma_*,\delta_*\}$ is one-sided repelling with  the interval $(\gamma_*,\delta_*)$ being the domain of immediate attraction of the fixed point $G_*$.
This fact makes the equilibrium $(G_*,I_*)$ locally attracting for the system (\ref{MS}). Outside the interval  $(\gamma_*,\delta_*)$ the structure of the map $\Phi$ can largely be preserved to be the same as in Theorem \ref{ThmA3}. This would guarantee the existence of a slowly oscillating periodic solution associated with the two-cycle $\{\alpha_*,\beta_*\}$. As in the case of Theorem \ref{ThmA3} the example can be easily generalized to produce any finite number of periodic solutions, while $(I_*,G_*)$ remains a locally attracting equilibrium. 
The phenomenon described by Theorem \ref{ThmA4} is demonstrated by an example of numerical solution of the system (\ref{MS-0}) in Subsection \ref{MPS}.

\bigskip

\subsection{Invariance, Persistence, and Global Asymptotic Stability}
\label{IPGAS}

Suppose that map $\Phi$ has a closed finite interval $L=[a,b]$ invariant in the general sense $\Phi(L)\subseteq L\subseteq\mathbb R_+$, and let interval $J$ be defined by $J=\tau_0f_1(L):=[c,d]$.
Consider the following subset $\mathbb X_L$ of the phase  space $\mathbb X$:
$$
\mathbb X_L=\{\psi=(\varphi,u)\in\mathbb X\,\vert\, u\in L, \varphi(s)\in J\; \forall s\in[-\tau,0]\}.
$$

It is easily seen, based of the properties of functions $F^{-1}$ and $H$, that the map $\Phi$ has an invariant interval $L$ such that for an arbitrary initial value $u_0\in\mathbb R_+$ its first iteration under $\Phi$, $u_1=\Phi(u_0)$ satisfies $u_1\in L$. Indeed, the interval $L$ can be defined as $L=[F^{-1}(H_{\infty}), F^{-1}(H_0)]$, where the finite interval $[H_0, H_{\infty}]$ is the image of the positive semi-axis $\mathbb R_+$ under the map $H$.
The values of $H_{\infty}$ and $H_0$ are given as $H_0=G_{in}+f_5(\tau_0a_0)),\; H_{\infty}=G_{in}+f_5(\tau_0a)$.
The corresponding interval $J$ is then defined as $J:=\tau_0 f_1(L)=[\tau_0f_1(F^{-1}(H_{\infty})),\tau_0f_1(F^{-1}(H_0)])]$.

The following statement describes the fact that the solutions of system (\ref{MS}) with initial functions in the set $\mathbb X_L$ remain within this set for all forward times $t\ge0$.

\begin{lem}\label{lem-inv}\,(Invariance)\,
Suppose that an initial function $\psi=(\varphi(s), u_0)$ is such that $\phi\in\mathbb X_{L}$, where $L$ is a closed interval invariant under map $\Phi$.
Then the corresponding solution $\mathbf x=\mathbf x(t,\psi)=(I(t),G(t))$ of system (\ref{MS}) satisfies $\mathbf x(t)\in\mathbb X_{L}$ for all $t\ge0$.
\end{lem}

Lemma \ref{lem-inv} shows that when the initial data for system (\ref{MS}) is such that $G(0)\in L$ and $I(s)\in J\; \forall s\in[-\tau,0],$ 
then the components $G$ and $I$ of the corresponding solution to system (\ref{MS}) satisfy the inclusions:
$$
G(t)\in L,\; I(t)\in J\quad\text{for all}\quad t\ge0.
$$
\begin{proof}
The proof of Lemma \ref{lem-inv} can be done by induction in time $t$ by using the cyclic structure of system (\ref{MS}). We provide its outline below.

Suppose that the initial function $\psi=(\phi(s), u_0)\in\mathbb X$ for system (\ref{MS}) is given such that $\phi(s)=I(s)\in J\;\forall s\in[-\tau,0]$ and $G(0)=u_0\in L$. 
Assume first that $G(t)\in L\;\forall t\in[0,T]$ for some $T>0$. Then also $I(t)\in J\; \forall t\in[0,T]$ is satisfied. 
Indeed, suppose $t_0\ge0$ is the first time moment of exit of the component $I$ from the interval $J$.
To be definite assume first that $I(t_0)=c$ and $I^\prime(t_0)<0$ and $I(t)<c\; \forall t\in(t_0,t_0+\varepsilon)$ for some $\varepsilon>0$.
Then $\tau_0f_1(G(t_0))\in J=[c,d]$ since $G(t_0)\in L=[a,b].$
Therefore, $\tau_0 I^\prime(t_0)=-c+\tau_0 f_1(G(t_0))\ge0$, a contradiction with $I^\prime(t_0)<0$.

In the case when $I(t_0)=c$ and $I^\prime(t_0)=0$ there exists a sequence $\{t_n\}$ of $t$-values such that $t_n\downarrow t_0$ and $I^\prime(t_n)<0, I(t_n)<c$.
This would imply that the derivative $I^\prime(t_n)=(1/\tau_0)[-I(t_n)+\tau_0 f_1(G(t_n))]>0$ is positive in a small right neighborhood of $t_0$, 
a contradiction with $t_0$ being the first point of exit from interval $J$.

Given $I(s)=\phi(s)\in J,\; \forall s\in[-\tau,0]$ and $G(0)=u_0\in L$ we shall show next that $G(t)\in L\; \forall t\in[0,\tau]$.
This is done in a way similar to the reasoning for the component $I$ above.
Assume $t_0\in[0,\tau]$ is the first point of exit of the component $G$ from the interval $L$.
To be specific let $G(t_0)=b$ and $G^\prime(t_0)>0$ holds. Using the monotone nature of functions $f_2$ and $f_4$ one sees that $f_2(G(t_0))+qG(t_0) f_4((I(t_0)))\le f_2(b)+q b f_4(b)$.
Therefore, $G^\prime(t_0)\le G_{in}+f_5(I(t_0-\tau))-f_2(b)-q b f_4(b)\le0,$ a contradiction with $G^\prime(t_0)>0$. 
The case when $G(t_0)=b, G^\prime(t_0)=0$ holds at the first point of exit from interval $L$ is treated similarly to the analogous case for $I(t)$ 
by selecting a sequence $t_n\downarrow t_0$ with $G(t_n)>b$ and $G^\prime(t_n)>0$.

The proof can now be completed by induction in $t$ with a step $\tau$.
Since $G(t)\in L\; \forall t\in[0,\tau]$ then also $I(t)\in J\; \forall t\in[0,\tau].$
These values of $G$ and $I$ are considered next as new initial data for the same solution to derive the inclusions 
$G(t)\in L, I(t)\in J,\; \forall t\in [\tau,2 \tau]$, and so on.
\end{proof}

From the proof of Lemma \ref{lem-inv} it is seen that for every initial data $\psi=(\phi(s),u_0)\in\mathbb X$ there exists a time moment $t=t_{\psi}$ such that the 
corresponding solution  $\mathbf x=\mathbf x(t,\psi)=(I(t),G(t))$ satisfies
\begin{equation}\label{incl_0}
I(t)\in J_0=\tau_0 f_1(\mathbb R_+)\quad\text{and}\quad G(t)\in L_0=\Phi(\mathbb R_+).
\end{equation} 
Indeed, if $I(t_0)\in J_0$ at some $t_0\ge0$ then $I(t)\in J_0\; \forall t\ge t_0,$ due to reasons in the first part of the proof of Lemma \ref{lem-inv}.
Likewise, $G(t)\in L_0\; \forall t\ge t_0$ if $G(t_0)\in L_0$ for some $t_0\ge0$.
Therefore, one has to consider the possibility that $I(t)\not\in J_0\; \forall t\ge0$ and $G(t)\not\in L_0\; \forall t\ge0$.
To be specific assume that $I(t)>\sup J_0$ and $G(t)>\sup L_0$ for all $t\ge0$ (other options are treated along the same line).
Then the respective equations of system (\ref{MS}) imply that $I^\prime(t)\le0$ and $G^\prime(t)\le0$ for all $t\ge0.$
Therefore, the finite limits $\lim_{t\to\infty} I(t)=I_\infty,\; \lim_{t\to\infty} G(t)=G_\infty$ exist.
By applying the limit to both equations of (\ref{MS}) along these components of the solution one sees that $(I_\infty, G_\infty)$ satisfies 
the equilibrium equations:
\[
f_1(G_\infty)=\frac{1}{\tau_0} I_\infty,\quad f_2(G_\infty)+q G_\infty f_4(I_\infty)=G_{in}+f_5(I_\infty).
\]
Therefore, $(I_\infty,G_\infty)$ is the only equilibrium of system (\ref{MS}), so that $I_\infty=I_*$ and $G_\infty=G_*$.
This is a contradiction with the inequalities $I_\infty\ge \sup J_0$ and $G_\infty\ge \sup L_0$, since $I_*$ and $G_*$ belong to the interior of the intervals $J_0$ and $L_0$, respectively.

The reasoning above leads to the following 

\begin{prop}\label{uni-persist} (Uniform Persistence I)
There are positive constants $0<m_I<M_I$ and $0<m_G<M_G$ such that for every initial data $\psi=(\varphi(s),u_0)\in\mathbb X$ there is 
a time moment $t=t(\psi)\ge0$ such that the corresponding solution $\mathbf x=\mathbf x(t,\psi)=(I(t),G(t))$ of system (\ref{MS}) satisfies
\[
m_I\le I(t)\le M_I\quad\text{and}\quad m_G\le G(t)\le M_G\quad \forall t\ge t_\psi.
\]
\end{prop}

Indeed, as it is seen from the above reasoning the values of the constants can be chosen as
\[
m_I:=\inf\{ \tau_0 f_1(\mathbb R_+)\},\; M_I:=\sup\{ \tau_0 f_1(\mathbb R_+)\},\; m_G:=\inf\{ \Phi(\mathbb R_+)\},\; M_G:=\sup\{ \Phi(\mathbb R_+)\}.
\]

We can now apply an inductive argument to the chain of reasoning preceding Proposition \ref{uni-persist}. Since $I(t)\in J_0$ and $G(t)\in L_0\; \forall t\ge t_0\ge0$ then 
$I(t)\in J_1=\tau_0 f_1(J_0)\subseteq \tau_0 f_1(\mathbb R_+)=\tau_0f_1(L_0)$ and $G(t)\in L_1=\Phi(L_0)\; \forall t\ge t_1\ge t_0.$
This is shown exactly  the same way as the inclusions (\ref{incl_0}). By the induction reasoning, there exists a sequence of $t$-values, $t_0\le t_1\le t_2\le \dots \le t_n\le t_{n+1}\le \dots,$ such that
\begin{equation}\label{imb_seq_Ln}
    I(t)\in \tau_0f_1(L_{n}):=J_{n+1}\quad{\text{and}}\quad G(t)\in L_{n+1}:=\Phi(L_n)\quad \forall t\ge t_{n+1}\ge t_n, \; n\in\mathbb N_0.
\end{equation}
The crucial role for the asymptotic behaviour of solutions $\mathbf x(t)=(I(t),G(t))$ is now played by the structure of the set $L_*=\cap_{n\ge0}\, L_n$.
Note that the imbedded sequence of intervals $L_0\supseteq L_1\supseteq L_2\supseteq \dots \supseteq L_n\supseteq L_{n+1}\supseteq\dots$, and the limiting set $L_*$ were constructed in Subsection \ref{PRs}.
The following two possibilities can only happen.

(I)\, The set $L_*=[\alpha_*,\beta_*]$ is a closed interval with non-empty interior.
Then points $\alpha_*<\beta_*$ form a cycle of period two under the map $\Phi$. 
It is the maximal cycle of period two for the map $\Phi$ in the sense that any other cycle of period two belongs to the open interval $(\alpha_*,\beta_*).$
Also, the cycle $\{\alpha_*,\beta_*\}$ is at least one-sided attracting (from above). The latter means that for every initial value $G_0\in(-\infty,\alpha_*)$ one has that 
$\Phi^{2n}(G_0)$ is an increasing sequence with $\lim_{n\to\infty} \Phi^{2n}(G_0)=\alpha_*$. Likewise, for every initial value $G_0\in(\beta_*,\infty)$ the sequence $\Phi^{2n}(G_0)$ is decreasing
with $\lim_{n\to\infty} \Phi^{2n}(G_0)=\beta_*$. Therefore, in this case, the persistence property of Proposition \ref{uni-persist} can be essentially improved.
Denote the interval $\tau_0 f_1([\alpha_*,\beta_*])=[c_*,d_*]$. The following property holds:
\begin{prop}\label{uni-persist2} 
     (Uniform Persistence II)\, For arbitrary initial data $\psi=(\varphi(s),u_0)\in\mathbb X$ the following holds for the corresponding solution $\mathbf x(t,\psi)=(I(t),G(t))$
     \[
     c_*\le\liminf_{t\to\infty} I(t)\le\limsup_{t\to\infty}I(t)\le d_*\quad\text{and}\quad \alpha_*\le\liminf_{t\to\infty} G(t)\le\limsup_{t\to\infty}G(t)\le \beta_*.
     \]
\end{prop}
The proof immediately follows from the property (\ref{imb_seq_Ln}). In fact,  more precise inequalities also hold under the assumptions of Proposition \ref{uni-persist2}:
\[
     c_*\le I(t)\le d_*\quad\text{and}\quad \alpha_*\le G(t) \le \beta_*\quad{\forall}\; t\ge t_*\ge0.
\]
A proof of the latter requires certain preliminaries and details which cannot be included in the paper due to their length.

(II)\, The set $L_*=[\alpha_*,\beta_*]$ is a single point.
This implies that $\alpha_*=\beta_*=G_*,$ and that the fixed point $G_*$ is globally attracting on $\mathbb R_+$ for the map $\Phi$.
In this case one has that the following global asymptotic stability property holds for system (\ref{MS}).

\begin{thm} (Global Asymptotic Stability, also Theorem~\ref{ThmA1}) 
\label{GAS}
Suppose that the unique fixed point $G_*$ of the interval map $\Phi$ is globally attracting: $\lim_{n\to\infty}\Phi^n(G)=G_*$ for every $G\in\mathbb R_+$.
Then the unique constant solution $(\tau_0f_1(G_*),G_*)$ of system (\ref{MS}) is globally asymptotically stable:
for arbitrary initial function  $\psi=(G(s),I_0)\in\mathbb X$  and every delay $\tau>0$ the following holds for the corresponding solution
$$
\lim_{t\to\infty}\,{\mathbf x}(t)=\lim_{t\to\infty}\,(I(t),G(t))=(\tau_0f_1(G_*),G_*)\,.
$$
\end{thm}
Again, the proof immediately follows from inclusions (\ref{imb_seq_Ln}).

{\bf Remark}.\,
    Note that a uniform persistence property of all solutions of system (\ref{MS}) is also proved in \cite{BenGou04-1}, see Proposition 2.4 there.
    However, our uniform persistence results, given by Propositions \ref{uni-persist} and \ref{uni-persist2}, provide explicit lower and upper bounds for the components $I$ and $G$ in terms of one-dimensional map $\Phi$
    (therefore, in terms of functions $f_1,f_2,f_4,f_5$ and parameters $\tau_0, q$). In fact, the bounds given by Proposition \ref{uni-persist2} are best possible in certain circumstances, e.g. when $\tau\to\infty$.
    They are given in terms of the maximal cycle of period two for the map $\Phi$.
    
    Paper \cite{BenGou04-1} also contains a condition for the global convergence to the equilibrium value  $G_*$ of the component $G(t)$ of system (\ref{MS}).
    It is given by Theorem 3.2 there, which requires that the following system for $x$ and $y$
    \begin{equation}\label{xy-cycle}
        G_{in}-f_2(x)-q x f_4(\tau_0 f_1(y))+f_5(\tau_0 f_1(y))=0,\quad G_{in}-f_2(y)-q y f_4(\tau_0 f_1(x))+f_5(\tau_0 f_1(x))=0
    \end{equation}
    has no solutions $x>0,y>0$. This is related to our more general and transparent condition of Theorem \ref{GAS}, about the global asymptotic stability in (\ref{MS-0}), which simply requires that the fixed point $G_*$ of the map $\Phi$ is globally attracting.
    If the later is satisfied then system (\ref{xy-cycle}) has no solutions $x>0,y>0$, since the existence of such a solution would mean that the pair $x,y$ forms a cycle of period two for the map $\Phi$, 
    contradicting the global attractivity of its fixed point $G_*$. In fact, it can be showed, with some additional effort, that under the assumptions imposed on system (\ref{MS}) the only fixed point $G_*$ 
    of map $\Phi$ is globally attracting if and only if system (\ref{xy-cycle}) has no positive solutions.

\subsection{Periodic Solutions}\label{PSs}


In this subsection we outline the algorithm how the existence of periodic solutions for system (\ref{MS}) can be derived. 
It follows the well established techniques of the ejective fixed point theory, see \cite{DieSvGSVLWal95} and \cite{HalSVL93} for general theoretical basics; we also use some related specific details from papers \cite{adH79a,BraHasIvaTro20,IvaDza20,IvaBLW04,IvaBLW19} to show the periodicity. 

The basic components for the existence of periodic solutions are:
\begin{itemize}
    \item[(1)] Construction of a cone of initial functions, and a translation operator along solutions  on it (Poincar{\'e} map), such that its fixed points give us slowly oscillating periodic solution. Some of these will have to be verified {\bf numerically};
    \item[(2)]
    The instability of the zero solution of the corresponding linearized system. This can be derived from the characteristic equation in terms of the existence of a pair of complex conjugate solutions with positive real part. Known results can be used here with proper harvesting and compilation, e.g. those in \cite{adH79a,BraHasIvaTro20,IvaDza20};
    \item[(3)]
    The compactness of the nonlinear map constructed in step (1) above. This is rather straightforward derivation based of the boundedness and smoothness properties of the nonlinear functions $f_1,f_2,f_4,f_5$ given in assumptions $(H1)-(H5)$;
    \item[(4)] Application of known results for the existence of periodic solutions for systems similar to (\ref{MS}). In particular, application of the well established ejective fixed point theory to our case;
\end{itemize}

The proof of existence of periodic solutions to system (\ref{MS}) (or equivalent system (\ref{MS-0})) uses well established theory of the ejective fixed point techniques applied to specially constructed maps on subsets of initial functions of the phase space. The subsets are usually cones of the initial functions generating the so-called slowly oscillating solutions. The related maps are appropriately constructed shifts along corresponding solutions. The general theory of such approach is described in e.g. \cite{DieSvGSVLWal95,HalSVL93}. In addition we shall use specific cases and results obtained in papers \cite{adH79a,BraHasIvaTro20,IvaDza20,IvaBLW04,IvaBLW19}.

\begin{defi}
(i) Given delay $\tau>0$ a continuous function $u(t): \mathbb R_+\to\mathbb R$ is called slowly oscillating (with respect to zero) if the distance between any two of its zeros is greater than $\tau$;\newline
(ii) A solution $(I(t),G(t))$ of system (\ref{MS}) is called slowly oscillating for $t\ge0$ if each of the functions $G(t)-G_*$ and $I(t)-I_*$ is slowly oscillating (with respect to zero in the sense of part (i)).
\end{defi}
In case when (ii) holds each of the components $G(t)$ and $I(t)$ is viewed as slowly oscillating function with respect to its constant component of the unique equilibrium $(G_*,I_*)$ of system (\ref{MS}).

We need a sufficient condition which guarantees the oscillatory nature of all solutions to system (\ref{MS}). We can use the corresponding result of paper \cite{BraHasIvaTro20}, see Theorem 1 there.

%
%
\begin{prop}\label{osc}
Suppose that nonlinearities $f_1,f_2,f_4, f_5$ are twice continuously differentialble on $\mathbb R$ and the characteristic equation (\ref{ChEq}) has no real solutions.
Then all solutions to system (\ref{MS}) oscillate about the positive equilibrium $(I_*,G_*)$.
\end{prop}
For the remainder of this subsection we shall assume that the conditions of Proposition \ref{osc} are satisfied.

{\bf Cone}. 
Consider the following set of initial functions $\mathbb K\subseteq \mathbb X$:
$$
\mathbb K=\{\psi=(\varphi(s),u)\in\mathbb X\; \vert\; u-G_*\ge 0,\, \varphi(s)-I_*\ge0,\;\text{and}\; \varphi(s)\exp\{1/\tau_0\, s\}\,\uparrow,\; s\in[-\tau,0) \}. 
$$
$\mathbb K$ is a cone on $\mathbb X$.

%
%
\begin{prop}\label{slow-osc}
Suppose that the characteristic equation (\ref{ChEq}) has no real solutions and the initial function $\psi=(\varphi(s),u)\in\mathbb K$ is such that $\varphi(s)\ge I_* \; \forall s\in[-\tau,0], \varphi(0)>I_*, u>G_*$. 
Then the corresponding solution $(I(t),G(t))$ of system (\ref{MS}) is slowly oscillating in the sense that each component $I(t)-I_*$ and $G(t)-G_*$ is slowly oscillating.
Moreover,
\begin{itemize}
    \item[(i)]
    The component $I(t)-I_*$ has a sequence of zeros $\{t_k\}$ such that $0<t_1<t_2<t_3<\cdots<t_k<t_{k+1}<\cdots $ and $t_{k+1}-t_k>\tau$ for all $k\in\mathbb N$.
    In addition, $I(t)-I_*<0$ for $t\in(t_{2k-1},t_{2k})$ and $I(t)-I_*>0$ for $t\in(t_{2k},t_{2k+1}),\, k\in\mathbb N$;
    \item[(ii)]
    The component $G(t)-G_*$ has a sequence of zeros $\{s_k\}$ such that $0<s_1<s_2<s_3<\cdots<s_k<s_{k+1}<\cdots $ and $s_{k+1}-s_k>\tau$ for all $k\in\mathbb N$.
    In addition, $G(t)-G_*<0$ for $t\in(s_{2k-1},s_{2k})$ and $G(t)-G_*>0$ for $t\in(s_{2k},s_{2k+1}),\, k\in\mathbb N$;
    \item[(iii)]
    The two sequences of zeros for $I-I_*$ and $G-G_*$ satisfy the following relationship: 
    $$
    s_1<t_1<s_2<t_2<s_3<t_3<\cdots< s_k<t_k<s_{k+1}<t_{k+1}<\cdots\, 
    $$
    with $s_{k+1}-t_k>\tau$ for all $k\in\mathbb N.$
\end{itemize}
\end{prop}
Main claims of Proposition \ref{slow-osc} are proved along the lines of similar propositions for other classes of equations; 
see e.g. \cite{HadTom77} for scalar equations, and \cite{IvaDza20,IvaBLW04,IvaBLW19} for systems.
We are still missing several details of a rigorous mathematical proof of this proposition; however, we have extensively verified it numerically for various choices of nonlinear functions $f_1,f_2, f_4, f_5$.

{\bf Mapping on Cone}. 
Proposition \ref{slow-osc} allows one to define a nonlinear map $\mathbb F$ on the cone $\mathbb K$ in the following way. 
Given initial function $\psi=(\phi(s),u)\in\mathbb K$ consider the corresponding solution $\mathbf x=(I(t),G(t)), t\ge0,$ to system (\ref{MS}). 
Given its second zero $s_2$ consider the first component $I(t)$ at time $s_2+1$ as an element $\phi_1(s)$ of the Banach space $C([-\tau,0],\mathbb R)$,
i.e $\phi_1(s):=I(s_2+1+s), s\in[-\tau,0]$. Then $\phi_1(s)>I_*\; \forall s\in(-\tau,0]$ and $u_1:=G(s_2+1)>G_*$, due to Proposition \ref{slow-osc}. Therefore, the mapping 
\begin{equation}
    \mathbb F: \psi= (\varphi(s),u)\mapsto (\varphi_1(s),u_1),
\end{equation}
maps cone $\mathbb K$ into itself. The mapping $\mathbb F$ is well defined for any $\psi\in\mathbb K$, different from the identical zero. 
For the trivial initial function $\psi\equiv(I_*,G_*)$ one defines $\mathbb F((I_*,G_*)):=\psi_1=(I_*,G_*)$, by the continuity of the map $\mathbb F$.

It is an easy observation that a nontrivial fixed point $\psi_0$ of the map $\mathbb F$, $\mathbb F(\psi_0)=\psi_0$, gives rise to a slowly periodic solution of system (\ref{MS}).
However, the map $\mathbb F$ always has the zero $\psi-(I_*,G_*)\equiv0$ as the trivial fixed point (which results in the identical zero solution to system (\ref{MS-0}) for $\forall\; t\ge0$). 
Therefore, one is interested in finding fixed points of map $\mathbb F$ which are different from the trivial zero one. 
This is done by application of the well developed theory of the ejective fixed point theory, which has been applied to various classes of functional differential equations elswhere.

{\bf Compactness and Boundedness}.
An important property required of map $\mathbb F$ in the ejective fixed point theory is its compactness and boundedness. 
It is a well known basic fact that that a shift operator along solutions of retarded differential delay equations is compact \cite{DieSvGSVLWal95,HalSVL93}.
The boundedness of $\mathbb F$ easily follows from the invariance property, Lemma \ref{lem-inv}.
One sees that for arbitrary $\psi=(\varphi,u)\in\mathbb K$ its first image under $\mathbb F, \psi_1=(\varphi_1(s), u_1)$ satisfies $\varphi_1(s)\in J_0$ and $u_1\in L_0$
(where the intervals $J_0,L_0$ are defined earlier).
Thus the set $\mathbb F(\mathbb K)$ is uniformly bounded from above and below.
Alternatively, one can also start with a bounded convex part $\mathbb K_0$ of cone $\mathbb K$, requiring that elements $\phi_0=(\varphi_0(s),u_0)\in\mathbb K_0$
satisfy $\varphi_0(s)\in J_0\;\forall s\in[-\tau,0]$ and $u_0\in L_0$.

{\bf Ejectivity}. 
The jectivity of map $\mathbb F$ can be determined in terms of a linear operator calculated on specific eigenvalues of the linearized system (\ref{LMS}) \cite{DieSvGSVLWal95,HalSVL93}.
It has a closed form in a general case \cite{IvaBLW19}. 
In lower dimensions of scalar equations or delay systems of two equations the property of ejectivity is eventually reduced to the existence of solutions of the characteristic equation (\ref{ChEq})
with positive real part and the imaginary part within the range $(0,\pi/\tau)$ \cite{adH79a,HadTom77,IvaDza20,IvaBLW04}.

We shall show next that for all sufficiently large delays $\tau$ the characteristic equation (\ref{ChEq}) has a pair of complex conjugate solutions $\alpha_0\pm i\beta_0$ with the positive real part $\alpha_0>0$ and the imaginary part $\beta_0$
satisfying $0<\beta_0<\pi/\tau$. This would imply the ejectivity of the above map $\mathbb F$.

It is more convenient to rewrite the characteristic equation (\ref{ChEq}) in an alternative form by rescaling the time $t=\tau\cdot s$ to get the normalized delay $\tau=1$ (see subsection \ref{lim_int_map}, system (\ref{MS-0})). One derives the following
\begin{equation}\label{ChEq_1}
    (\varepsilon\lambda+\mu_1)(\varepsilon\lambda+\mu_2)+b+a\,\exp\{-\lambda\}=0,
\end{equation}
where $\varepsilon=1/\tau>0$ is a small parameter when $\tau>0$ is large enough.
By setting $\varepsilon=0$ one gets the equation $\mu_1\mu_2+b+a\,\exp\{-\lambda\}=0$, which has a pair of complex conjugate solutions $\lambda=\alpha_0\pm i\pi$, where $\alpha_0=\ln[a/(\mu_1\mu_2+b)]>0$.
Consider now the characteristic equation (\ref{ChEq_1}) for small $\varepsilon>0$. 
By Rouch\'{e}'s Theorem  it has a pair of complex conjugate solutions $\lambda_\varepsilon=\mu(\varepsilon)\pm i\nu(\varepsilon)$ 
such that $\mu(\varepsilon)$ is close to $\alpha_0>0$ and $\nu(\varepsilon)$ is close to $\pi$.
We shall show that $\nu(\varepsilon)<\pi$ for all sufficiently small $\varepsilon>0$.
One rewrites the characteristic equation (\ref{ChEq_1}) for the solution $\lambda_\varepsilon$ in the form
$$
(\varepsilon\mu+\mu_1+\varepsilon\nu i) (\varepsilon\mu+\mu_2+\varepsilon\nu i)+b+a\,\exp\{-\mu\}(\cos\nu-i\sin\nu)=0.
$$
and considers its imaginary part:
$$
\varepsilon\nu(2\varepsilon\mu+\mu_1+\mu_2)-a\,\exp\{-\mu\}\sin\nu=0.
$$
By differentiating the last equation with respect to $\varepsilon$ and setting $\varepsilon=0$ one finds
$$
\nu^\prime(0)=-\frac{\pi(\mu_1+\mu_2)}{\mu_1+\mu_2+b}<0,
$$
which proves that $\nu(\varepsilon)<\pi$ for all sufficiently small $\varepsilon>0$, since $\nu(0)=\pi$.

\subsection{Multiple Periodic Solutions}\label{MPS}

We will demonstrate numerically the existence of multiple periodic solutions using system~(\ref{MS-0}). We start with linear functions $F_1$, $F_2$, and $F_4$, which contain a constant function $f_4$. The only non-linear function is then $F_5(x) = F(x)$, which is monotonically decreasing (see system (\ref{MS-s}) below). The two-dimensional system of this type is simply looking and close in a sense to a single scalar differential delay equation where the non-uniqueness of slowly periodic solutions is known by several publications \cite{IvaSha91,Nus79b}. Having derived multiple periodic solutions  for system (\ref{MS-s}) we will perturb it by the inverse tangent function to produce a system of type (\ref{MS-in}) which will have the same two periodic solutions, however, they are slightly perturbed compared with those in system (\ref{MS-s}).

With the first step, system~(\ref{MS-0}) becomes
\begin{eqnarray}\label{MS-s}
            x^\prime(t) &=& -\frac{1}{\tau_0} x(t) + a_1 y(t)  \\
  \nonumber y^\prime(t) &=& -a_2 y(t) - a_4 x(t) + F(x(t-\tau)),
\end{eqnarray}
where $a_i > 0$ for $i=1,2,4$. Here, we will consider  two appropriate choices of the monotonically decreasing $F(x)$ designed  as follows:
\begin{equation}
\label{Fa}
    F(x) =
    \begin{cases}
    f(x), &  x \in [0, M] \\
    -x, & x \in [M, \frac{\pi}{2}] \\
    -\frac{\pi}{2} - A \arctan [k(x-\frac{\pi}{2})], & x > \frac{\pi}{2} \\
    -F(-x) & x < 0.
    \end{cases}
\end{equation}
where, $A > 0$ and $k > 1$ are positive arbitrary constants. For this definition,  $M$ is the solution of the equation $\arctan(x)=x$. 
Therefore, the function $F$ is continuous (but not $C^1$) and odd by construction.

\begin{figure}
    \includegraphics{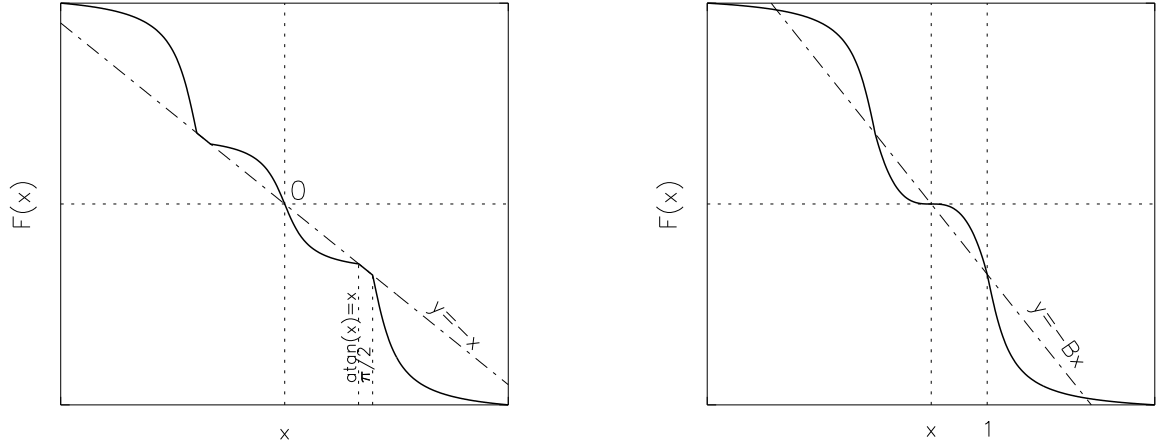}
    \caption{Functions $F(x)$ used to demonstrate multiple solutions behaviour: left plot - $F(x)$ defined by Equation~(\ref{Fa}), right plot - $F(x)$ defined by Equation~(\ref{Fb}).}
    \label{fig:f5_functions}
\end{figure}

Another choice of function $F$ is as follows:
\begin{equation}
\label{Fb}
    F(x) =
    \begin{cases}
    -Bx^{2n+1}, &  |x| \le 1,~B > 0 \\
    -B - A \arctan [k(x-1)], & x \ge 1 \\
    -F(-x) & x \le -1.
    \end{cases}
\end{equation}
The function plots are shown in Fig.~\ref{fig:f5_functions}. These choices allow us to demonstrate the presence of multiple periodic solutions of different types.  The first choice (Equation~(\ref{Fa})) leads to two different periodic solutions, while the second choice (Equation~(\ref{Fb})) to an attracting equilibrium and oscillating periodic solutions of  system~(\ref{MS-0}), as demonstrated by numerical solutions in Section~\ref{numerical}. 

We verify Theorem~\ref{ThmA1} numerically in Section~\ref{numerical}. We make a small modification to  system~(\ref{MS-s}) so that it can be viewed as the original system of the form (\ref{MS-0}). 

We change  function $f_4$ from a constant to a monotonically increasing function with $0<d = f_4(0) < \lim_{x \rightarrow \infty} f_4(x) = e > d$. For such $f_4(u)$ one can choose:

\begin{equation}
    f_4(u) = \epsilon (A + B \arctan (u)),~u \in\mathbb R,
    \label{f4_eq}
\end{equation}
where $A$, $B$ and $\epsilon > 0$ are constants, and $A > \frac{\pi}{2}B$.

We consider an intermediate system (between~(\ref{MS-0}) and (\ref{MS-s})), as follows:
\begin{eqnarray}\label{MS-in}
            x^\prime(t) &=& -\frac{1}{\tau_0} x(t) + a_1 y(t)  \\
  \nonumber y^\prime(t) &=& -a_2 y(t) - f_4(x(t)) y(t) -\delta B \arctan(x(t)) + F_5(x(t-\tau)),
\end{eqnarray}
where $\delta$ is a constant comparable to $\epsilon$ and $F_5$ is chosen as $F(x)$ from equations~(\ref{Fa}) and (\ref{Fb}).

By replacing the constants $a_1$ and $a_2$ in equation~(\ref{MS-in}) with  non-linear piecewise continuous functions, linearly proportional to the argument in its range within the span of the periodic solutions and equal to constants outside this range, and by replacing $F_5(x)$ outside the range of $x(t)$ by a symmetric smooth nonlinearity with a finite limit $\lim_{x\to\infty} F(x)=-F_\infty=-\lim_{x\to-\infty} F(x)$ for an appropriate $F_\infty>0$, the system (\ref{MS-in}) is converted back to the original form (\ref{MS-0}). 
   
    {\bf Remark}.\,
    The existence of any number of stable slowly oscillating  periodic solutions can be achieved in two different ways.
    
    (i) An analogous construction to that of function $F(x)$ given by (\ref{Fa}) can be continued on the interval beyond the amplitude of the second large periodic solution. Indeed, given $A>0$ and $k>1$ such that the second periodic solution exists, one finds the unique value $M_1>\pi/2$ such that it solves the equation $\pi/2+\arctan [k(x-\pi/2)]=x$. Then one defines function $\tilde{F}, x\ge0,$ such that $\tilde{F}\equiv F(x)$ for $x\in[0,M_1]$ and $\tilde{F}=-M_1-A_1\,\arctan[k_1(x-M_1)]$ for $x\ge M_1$, and $\tilde{F}(x)=-\tilde{F}(-x)$ for $x<0$. Exactly as with $F(x)$ given by (\ref{Fa}) it can be showed for the modified $\tilde{F}(x)$ that there exists $k_1^0$ large enough such that system  (\ref{MS-in}) has three slowly oscillating periodic solutions, with the amplitude of the largest one greater than $M_1$. This procedure of additional modification of $F(x)$ in (\ref{Fa}) can be continued step-by-step further so that one can obtain any finite number of stable slowly oscillating periodic solutions. If the procedure is applied to function $F(x)$ given by (\ref{Fb}) then one derives any number of stable periodic solutions together with the locally stable equilibrium. 
    
    (ii) It is known that the existence of stable (hyperbolic in general) slowly oscillating periodic solutions persists under small continuous perturbations of the non-linear right hand side (functions $F, F_5, f_4, \arctan(\cdot)$ and constants $a_1,a_2,a_4, \tau_0$ for  systems (\ref{MS-s}) and (\ref{MS-in})) \cite{BLW1,BLW2}. Therefore, if the nonlinearity $F$ in (\ref{Fa}) is replaced, in a sufficiently small neighborhood $|x|<\delta$ of $x=0$, by an arbitrary and small function $\tilde{F}$, and $F$ remains the same outside the small vicinity, for $|x|\ge\delta$, then the two stable slowly oscillating periodic solutions will persist, having changed only a little. The replacement of $F(x)$ for $|x|<\delta$ can be done in such a way that the resulting function $\tilde{F}(x)$ is monotone decreasing there (therefore, it is monotone decreasing for all $x\in\mathbb{R}$). We now consider function $F(x)$ by (\ref{Fa}) on the interval $[-M_1,M_1]$ where $M_1$ is defined above in part (i). Rescale it next to the interval $[-\delta,\delta]$ by $\tilde{F(x)}=(\delta/M_1)F(\frac{M_1}{\delta} x)$. We now use the above $\tilde{F}$ to replace the original $F$ in the delta neighborhood of $x=0$. The resulting nonlinearity is now such that the corresponding system (\ref{MS-in}) has four stable slowly oscillating periodic solutions: two are the perturbed original periodic solutions, and the other two are small scaled original periodic solutions placed in the $\delta$-neighborhood of $x=0$. This procedure can be repeated any finite number of times.   
    

\section{Numerical Analysis}
\label{numerical}

Analytical investigation of systems of delay-differential equations and, in particular, system~(\ref{MS}), with biologically-inspired functions and experimentally measured parameters, is usually very difficult or impossible. Therefore numerical methods have to be employed to study the details of behaviours of the glucose-insulin regulation models \cite{engelborgs2001}. Li, Kuang and Mason \cite{LiKuaMas06} performed numerical analysis of a two-delay glucose-insulin regulation system to analyse the dependence of bifurcations in the system on delays. This model utilised functions $f_1-f_5$ in their exponential forms with experimentally determined constants. In papers \cite{HuaEasAng15,HuaBriAng17}  numerical analyses are performed on a similar system with more complex Hill functions, allowing for more realistic modelling of the physiological mechanisms of glucose-insulin regulation. They also studied the sensitivity of the solutions to the values of Hill parameters used and performed simulations, which represented glucose-insulin regulation disorders, namely both Type 1 and Type 2 diabetes. Here we also use numerical analysis to further clarify some of the analytical results, obtained in the previous sections. 

There are two main points we aim at to demonstrate numerically. First, we demonstrate usability of equation~(\ref{DE}) in diagnostics of the solution behaviour of system~(\ref{MS}). Then we revisit the statement on relative insignificance of the actual forms of functions $f_1-f_5$ \cite{keenerandsnayd} in comparison to their shapes.

\subsection{Numerical Methods}\label{NMs}

To confirm the results obtained in the previous sections, we produce numerical solutions for systems~(\ref{MS}), (\ref{MS-0}), (\ref{MS-in}) and equation~(\ref{DE}). Furthermore, some of the theoretical concepts and results obtained in Section~\ref{main_results} cannot be proven analytically, therefore we use numerical methods to verify their validity.

\begin{figure*}
\includegraphics{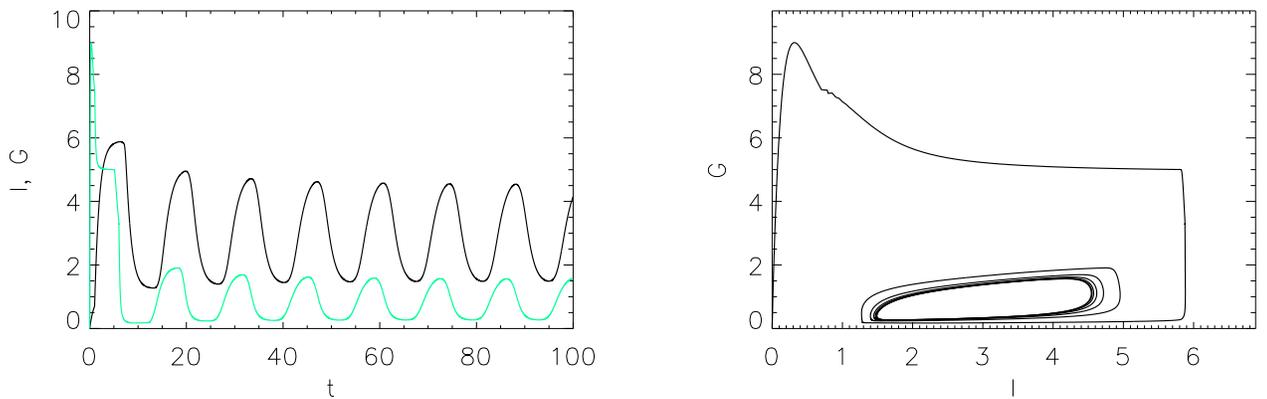}
\caption{An example of periodic solution to the system (\ref{MS}). The time evolution of $I$ (black) and $G$ (green) is shown in the left panel. Right panel shows the corresponding phase portrait for the system, plotted for a larger time interval $0 < t < 200$.}
\label{fig_oscil_soln}
\end{figure*}

\begin{figure*}
\includegraphics{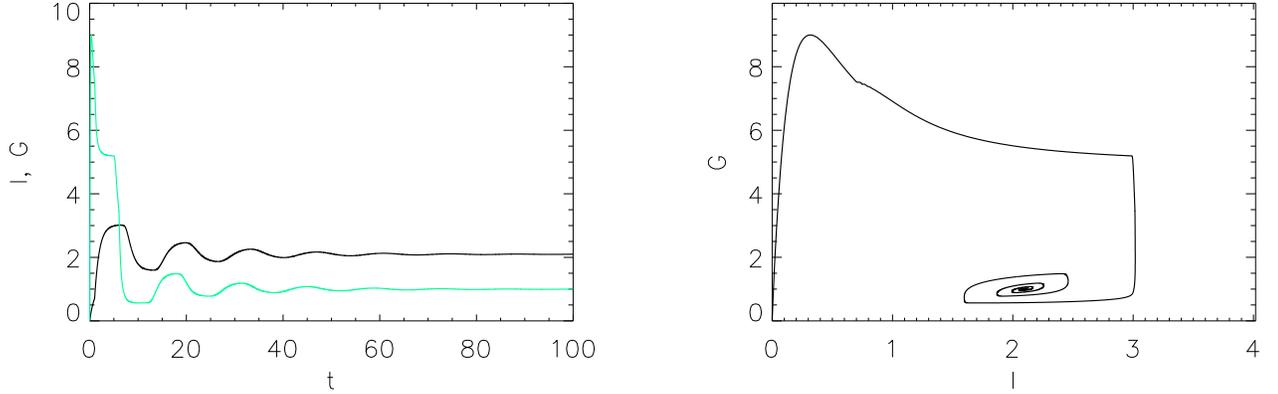}
\caption{Same as in Fig.~\ref{fig_oscil_soln}, but for an asymptotically stable solution to the system (\ref{MS}).}
\label{fig_stable_soln}
\end{figure*}

The initial value problem to system~(\ref{MS}) is solved by using a 4-th order Runge-Kutta-Fehlberg method with an adaptive time step. The solution examples and their corresponding phase portraits are shown in Figures~\ref{fig_oscil_soln} and \ref{fig_stable_soln}, which represent a periodic and an asymptotically stable solutions, respectively.

\begin{figure}
    \includegraphics{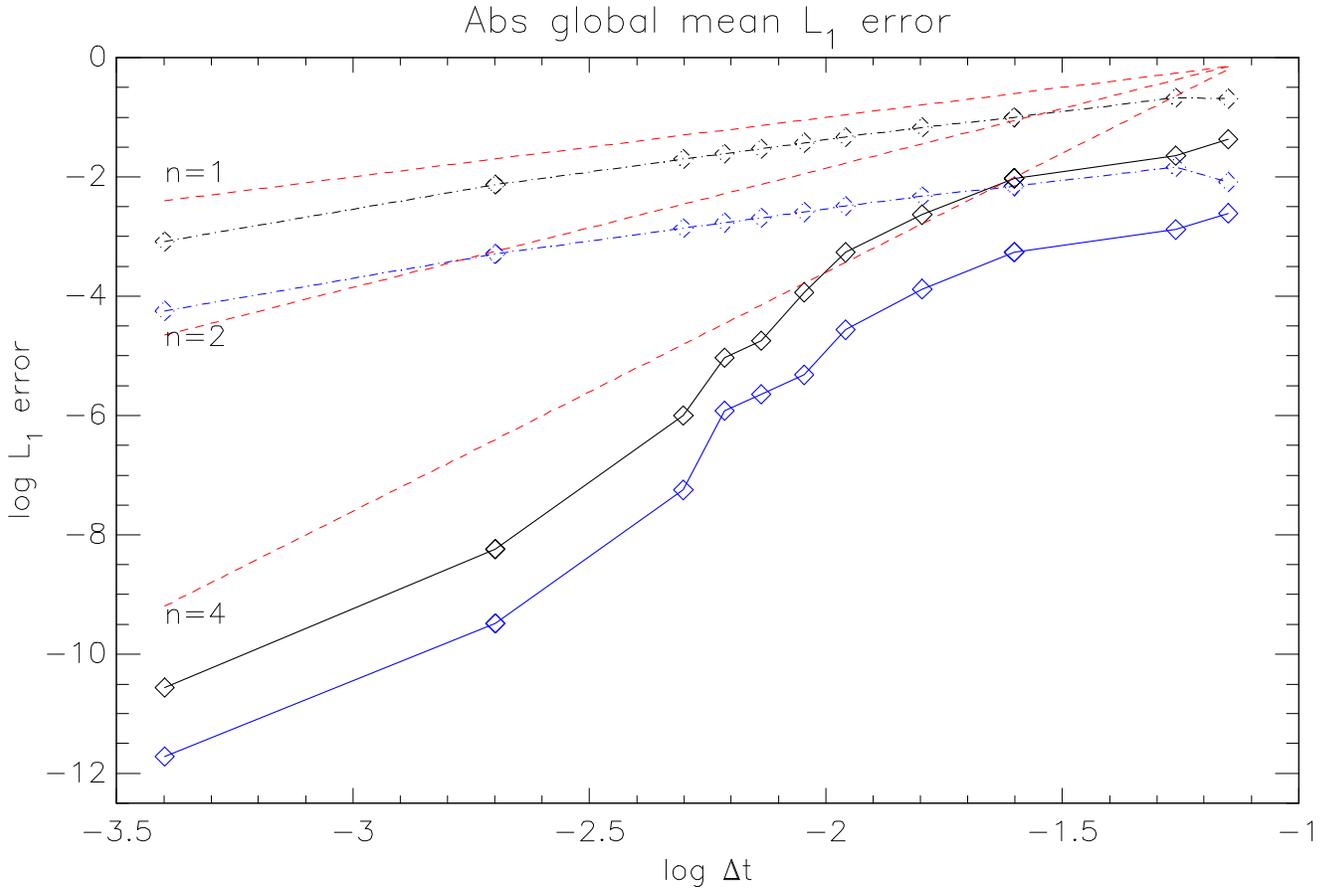}
    \caption{Dependence of the absolute global mean $L_1$ error on the time step for the employed numerical scheme. To demonstrate the precision order, the red dashed lines correspond to the power laws with the provided indices. Blue and black dash-dotted curves show Euler integration of the system with 4-th order Lagrange-interpolated delay term for $I$ and $G$, respectively. The solid curves show the 4-th order Runge-Kutta integration.}
    \label{fig:convergence}
\end{figure}

The delay term in the system is interpolated using Lagrange polynomials in their barycentric form \cite{BerTre2004}. This method demonstrates 4-th order self-convergence for sufficiently small time steps for both periodic and asymptotically stable solutions and a wide range of delays (see Figure~\ref{fig:convergence}).

To demonstrate applicability of the limiting interval map analysis, described in Subsection~\ref{lim_int_map}, we numerically solve Equation~(\ref{DE}). The solution of the (implicit with respect to $G(s)$) difference equation~(\ref{DE}) is preferential for numerical treatment as it does not require numerically inverting a function on an arbitrary range of its argument, despite equation~(\ref{IM}) being mathematically simpler and providing an explicit solution for $G(s)$.

Since the functions $f_1$-$f_5$ are monotone, numerical solution of the difference equation~(\ref{DE}) for $G(s)$ does not represent difficulties, and a simplest bisection method has been implemented. To distinguish numerically the solution types is also straightforward, as the period of the solution (if such period exists) is always 2 by construction. A solution is considered periodic for a large integer $s$ if $|G(s+2)-G(s)| < \epsilon$ and $|G(s+1)-G(s)| > \epsilon$, where $\epsilon=10^{-3}$ is a constant, which determines the precision.

\begin{figure*}
\includegraphics{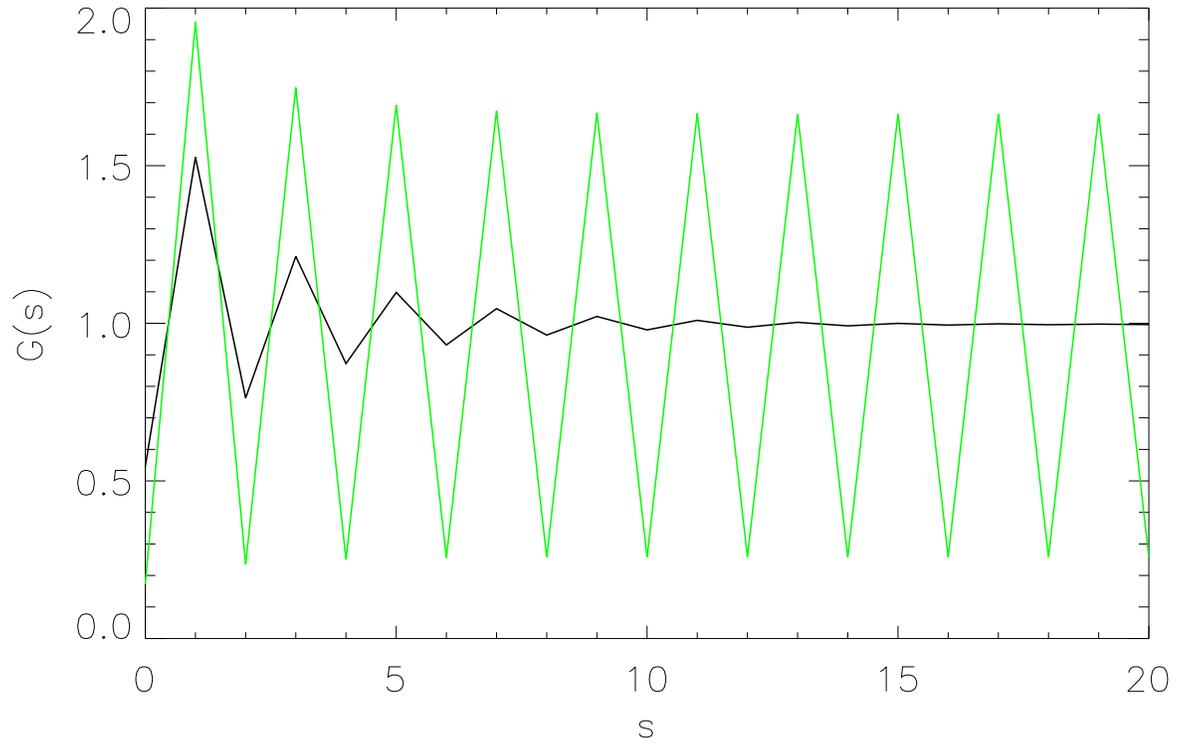}
\caption{Examples of periodic (green) and asymptotically stable (black) solutions of Equation~(\ref{DE}).}
\label{fig_DE_solns}
\end{figure*}

The solution of equation~(\ref{DE}) either exhibits an asymptotic stability, which corresponds to the asymptotically stable regime for any delay $\tau$ in system (\ref{MS}), or an oscillatory function with a period 2. The latter case corresponds to the periodic solution of system~(\ref{MS}), which exists for the delay $\tau$ greater than some critical value $\tau_c$, determined numerically from the full solution of system~(\ref{MS}) given a set of its parameters. If $\tau < \tau_c$, the system shows a stable equilibrium solution. Examples of  solutions to equation~(\ref{DE}) are shown in Fig.~\ref{fig_DE_solns}.

Other advantages of using equation~(\ref{DE}) in comparison to the original system~(\ref{MS}) are that it does not explicitly contain the delay value, neither does it require a priori knowledge of the oscillation period (if present) and the solution derivatives. It is, therefore, beneficial to numerically analyse the system's behaviour using this equation.

\begin{figure*}
\includegraphics{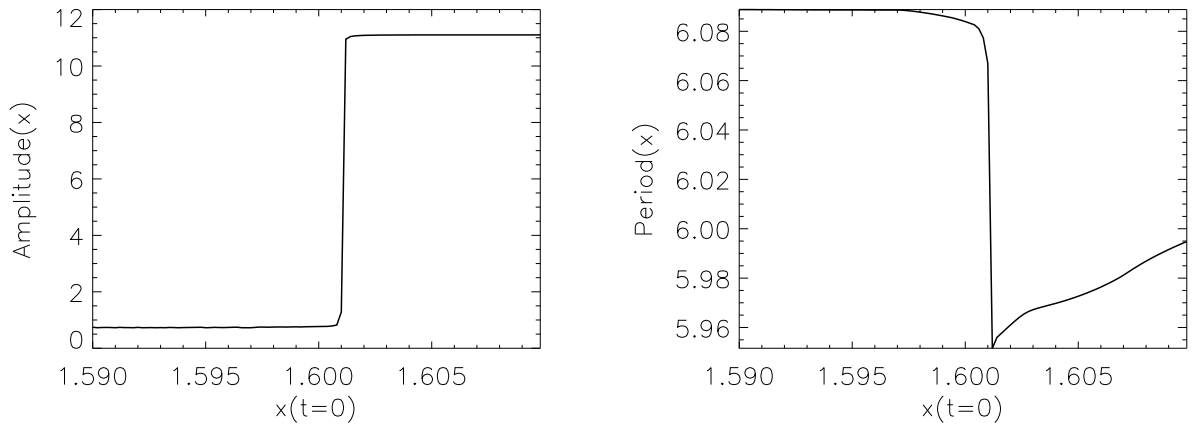}
\caption{Demonstration of multiple solutions of the system~(\ref{MS-in}) with $F$ as defined by equation~(\ref{Fa}). Left panel: dependence of the amplitude of the solution on the initial value $x(t=0)=y(t=0)$. Right panel: dependence of the period of the solution on the initial value $x(t=0)=y(t=0)$.}
\label{initial_condition_test}
\end{figure*}

\begin{figure}
\includegraphics{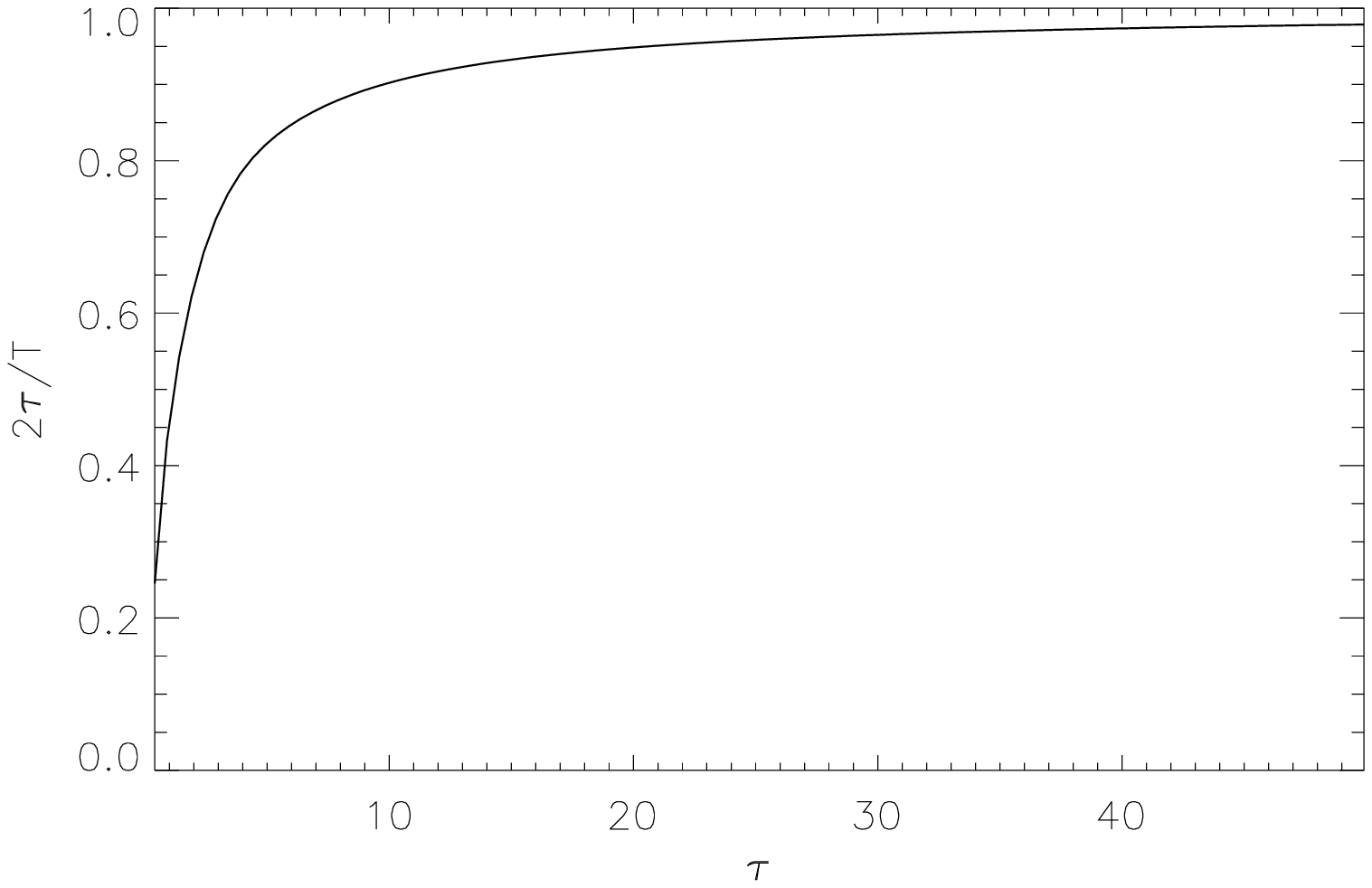}
\caption{Dependence of the ratio of the delay $\tau$ to half-period of the solution of system~(\ref{MS-in}) with $F_5$ as defined by equation~(\ref{Fa}) on the delay $\tau$, which confirms the slowly oscillating solution property.}
\label{tau_vs_T}
\end{figure}

\subsection{Numerical Demonstration of Multiple Periodic Solutions and Slow Oscillations}\label{NDMPSSO}
\begin{figure}
\begin{subfigure}{.5\textwidth}
  \centering
  \includegraphics[width=.8\linewidth]{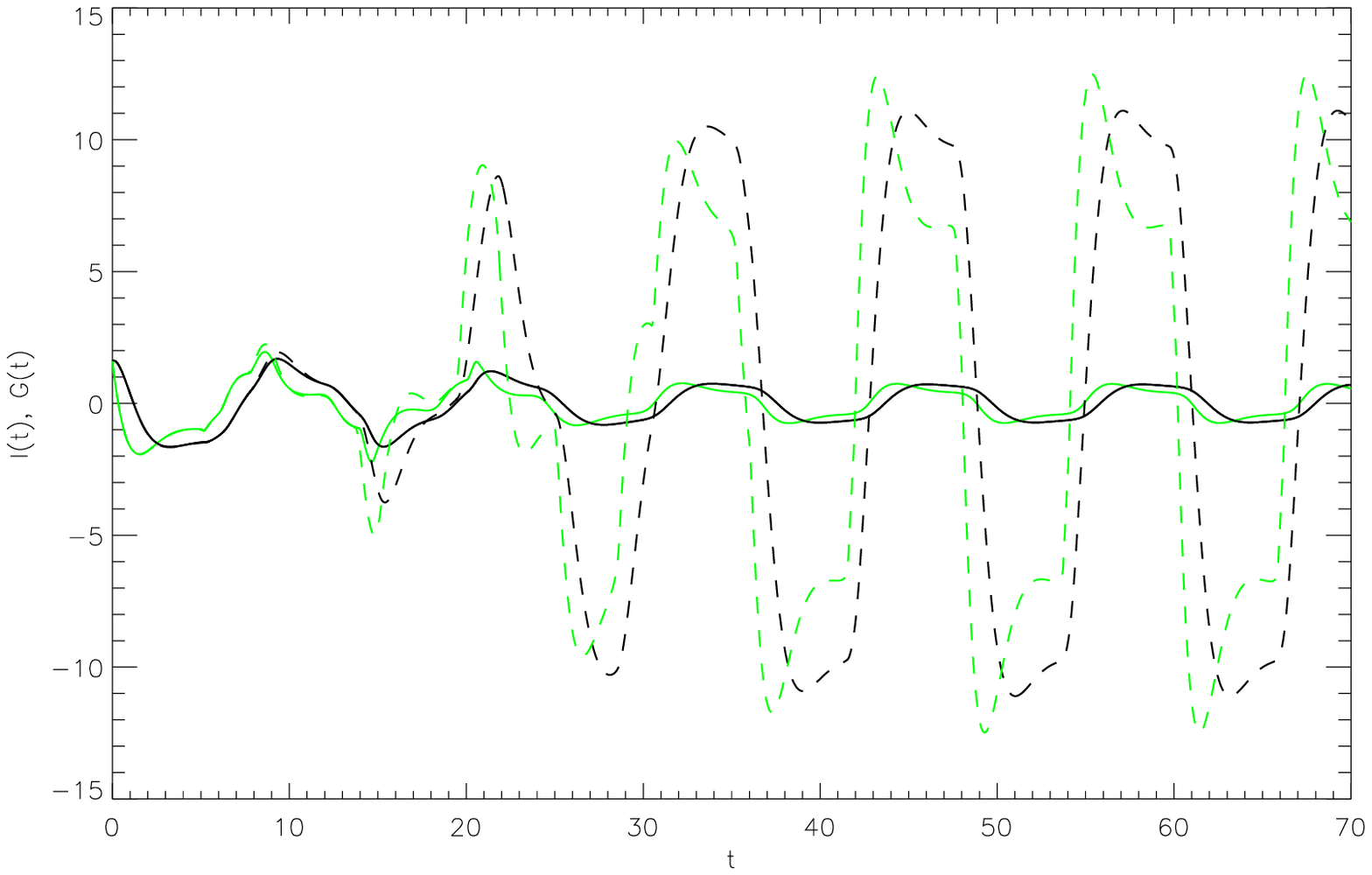}
  \label{fig_multiple_oscil}
\end{subfigure}%
\begin{subfigure}{.5\textwidth}
  \centering
  \includegraphics[width=.8\linewidth]{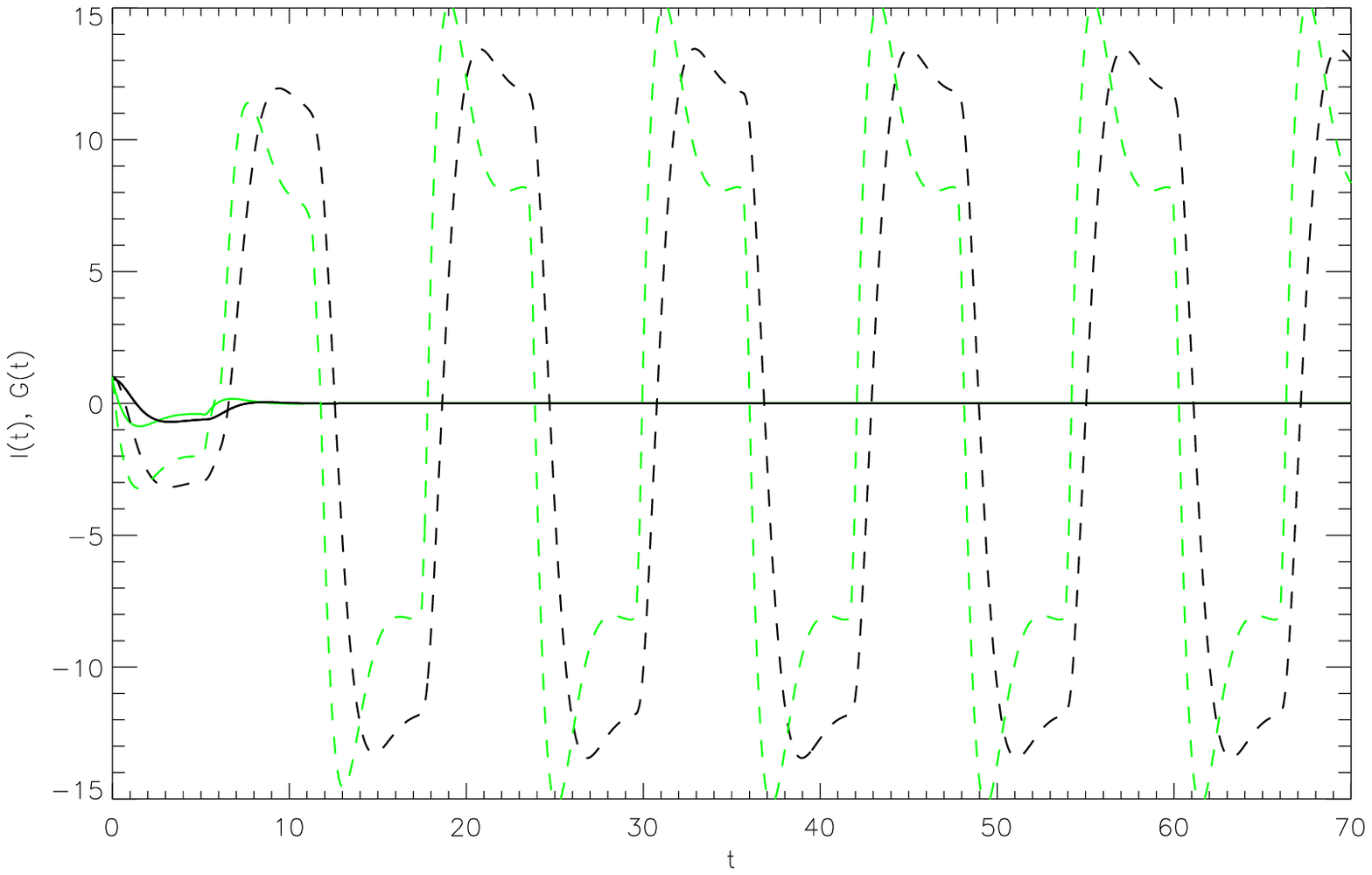}
  \label{fig_multiple_os_nos}
\end{subfigure}
\begin{subfigure}{.5\textwidth}
  \centering
  \includegraphics[width=.8\linewidth]{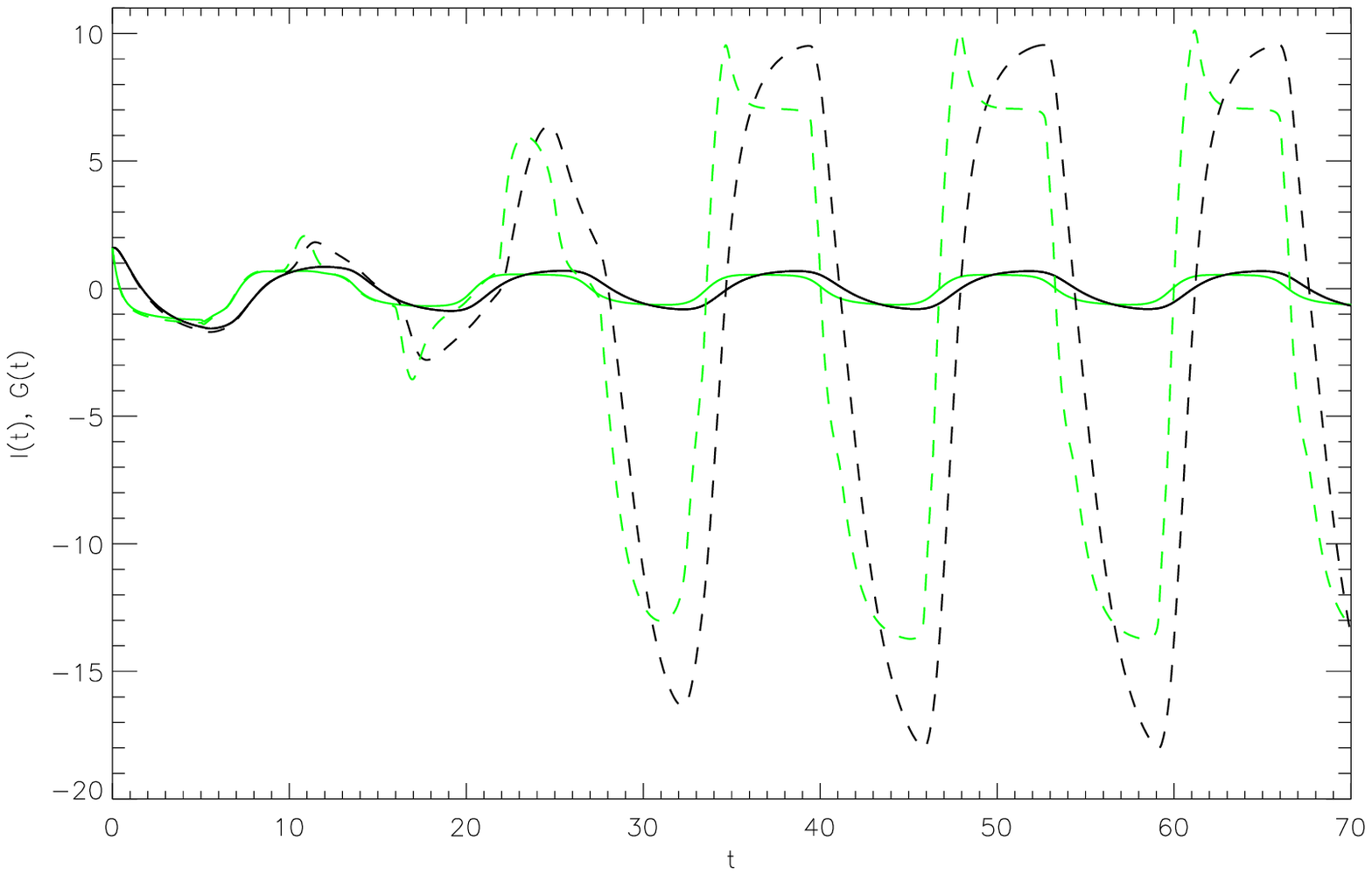}
  \label{fig_multiple_oscil_in}
\end{subfigure}
\begin{subfigure}{.5\textwidth}
  \centering
  \includegraphics[width=.8\linewidth]{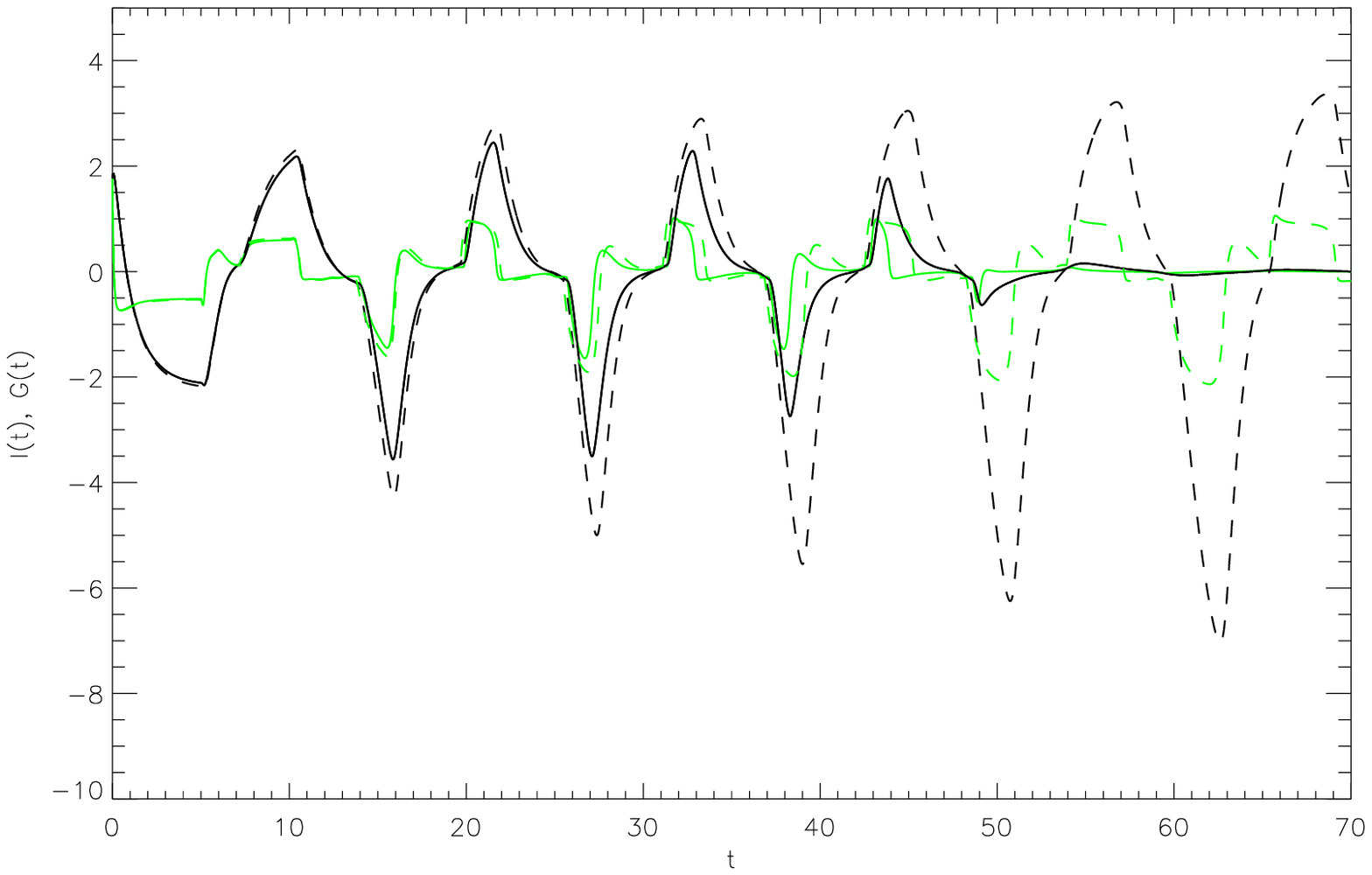}
  \label{fig_multiple_os_nos_oscil_in}
\end{subfigure}
\caption{Numerical demonstration of multiple periodic solutions of the two-dimensional system. Top and bottom rows show examples of multiple solutions for systems~(\ref{MS-s}) and (\ref{MS-in}), respectively. Left and right columns demonstrate multiple periodic and periodic/attractive equilibrium solutions, produced by functions $F$, defined by equations~(\ref{Fa}) and (\ref{Fb}). Green curves correspond to $G(t)$ and black curves correspond to $I(t)$. The delay $\tau=5$ is used for all solutions in the plot.}
\label{fig_multiple_solutions}
\end{figure}

To further verify Theorem~\ref{ThmA1}, we solve the system~(\ref{MS-in}) numerically. The piecewise functions $f_1$ and $f_2$, constructed as described above, and $f_4$ as in equation~(\ref{f4_eq}), are used in the calculation. Two different cases are considered for $F(x)$, as given in equations~(\ref{Fa}) and (\ref{Fb}), leading to different solution types. In Figure~\ref{fig_multiple_solutions}, examples of the solutions are shown. Transition between the different solutions of system~(\ref{MS-in}) occurs in a very narrow range of the initial conditions $x(t=0)=y(t=0)$. Figure~\ref{initial_condition_test} demonstrates the solution amplitude (left panel) and the solution period (right panel) for $x(t=0)=y(t=0)=[1.59,1.61]$. This figure also shows that there is a small effect (2\%) of the initial condition on the period of oscillations, with the transition occurring at the same value as the transition between the amplitudes of the solutions.

On the other hand, the time delay $\tau$ determines the period $T$ of oscillations.  This is illustrated in Figure~\ref{tau_vs_T}, where the dependence of ratio of the delay $\tau$ to the oscillation half-period $2\tau/T$ vs $\tau$ is shown. For all reasonable from the practical point of view values of $\tau$, $T \gtrsim 2\tau$. This confirms the existence of slow oscillations for this system. Recall that an oscillation is considered to be a slow oscillation if its half-period is greater than the delay $\tau$.
As system~(\ref{MS-in}) mimics the behaviour of the original system~(\ref{MS}), this shows that the time delay to a great extent determines the period of the oscillations and slow oscillations occur.



\section{Conclusion}\label{Conclusion}
In this paper we performed a further analytical and numerical study of the Sturis-Bennett-Gourley model describing the glucose-insulin regulation system in humans. The model is given by a nonlinear two-dimensional system of delay-differential equations with a single delay.  One of the principal goals of the paper was to demonstrate the applicability of the limiting interval maps approach to provide information on the system's asymptotic behaviour and to show the existence of slowly oscillating periodic solutions when the unique equilibrium is unstable.  

The model was introduced in papers \cite{BenGou04-1, BenGou04-2, BenGou04-3}; it includes one delay - namely the delay between plasma insulin production and its effect on hepatic glucose production. This model was selected to demonstrate the power of limiting interval maps method.

The method not only reproduced some of the results obtained in \cite{BenGou04-1, BenGou04-2, BenGou04-3}, but also showed the rich behaviour of the system with a  choice of physiological functions $f_i$, with specific attention on $f_5$. We investigated the behaviour of the system with $f_5(u)$ chosen as monotonically decreasing function in $u\in\mathbb R_+$ and showed that this specific choice leads to multiple oscillating periodic solutions or stable solutions converging to the equilibrium. We have demonstrated that depending on the appropriate choices of functions $f_1-f_5$ (which still satisfy all the conditions H1-H5), the non-uniqueness of the periodic solutions and their coexistence with the stable equilibrium can be achieved.

We would like to further notice  that our global asymptotic stability result Theorem \ref{ThmA1} can likely be deduced from considerations in paper \cite{BenGou04-1}. In particular, Theorem 3.2 there provides a sufficient condition for the global attractivity of the component $G(t)$ for a single integro-differential equation which is a truncated version of system (\ref{MS}). The component $G$ being attracted by its equilibrium value $G_*$ immediately implies, via an integral representation of the first equation of system (\ref{MS}), that the other component $I$ is attracted by its respective constant value $I_*$. One principal issue with the relevant considerations in \cite{BenGou04-1} is that integro-differential equation (3.1) is not an exact reduction but an approximation to system (\ref{MS}), which is achieved by dropping exponentially small perturbation terms.

In all our numerical simulations the eventual periodic solutions appear to be of the so-called "sinusoidal type". This means that they have a single maximum and a single minimum values and are monotone in between on their period. This shape of periodic solutions is rigorously proved in \cite{M-PSel96b} for the so-called unidirectional systems. System (\ref{MS}) is not of the unidirectional type, so this result cannot be directly extended to our case.

The model with single delay has been succeeded by a number of more sophisticated models with two delays \cite{LiKuaMas06, LiKua07, LiWang2012, HuaEasAng15}, which involve control loops containing muscle \cite{kissler2014} or effect of diabetes type I or II \cite{HuaBriAng17}. However, as one of the delays is always significantly larger than the other, a system with one delay can be a very good approximation to those with two delays. 

Constructing the nonlinear maps, we have found a difference equation, which represents the dynamics of the system in the large delay limit, which has the potential for  diagnostics of  the solution types without the need to solve the full system of differential equations with one delay.

The paper shows the elegance and efficiency of the approach via limiting interval maps in solving systems of differential equations with one delay. Furthermore, using this method, we revealed the existence of multiple slowly oscillating  periodic solutions, their coexistence with the stably equilibrium, or the global asymptotic stability of the unique equilibrium.

Thus, the paper shows the potential of this method for solving complex problems in mathematical physiology and is generally applicable for the systems of nonlinear differential equations with a single delay. 

\subsection{Author Contributions}
MA proposed the idea for the investigation. AI developed the theoretical aspects of the paper. SS produced the numerical results. GB contributed to the numerical aspects of the paper. All authors contributed to writing up the manuscript.

\subsection*{Acknowledgement}
This work was initiated during AI's visit to Deakin University, Burwood Campus, in December 2017. He would like to express his appreciation of the accommodation  and support from the School of Information Technology, Faculty of Science, and of the hospitality and collegiality from staff and his coauthors. This research was undertaken with the assistance of resources and services from the National Computational Infrastructure (NCI), which is supported by the Australian Government. MA thanks Newton Advanced Fellowship (UK Royal Society) / Academy of Medical Sciences UK for partial funding to develop this research.

\vspace{5mm}



\begin{thebibliography}{99}

\bibitem{adH79a}
an der Heiden, U. Periodic solutions of a nonlinear second order differential equation  with delay.  J. Math. Anal. Appl. 1979, 70, 599-609.

\bibitem{BelCoo63}
Bellman, R. and K.L. Cooke, K.L. Differential-Difference Equations. Academic Press, New York/London, 1963.

\bibitem{Bendixson} Bendixson, I. Sur les courbes d{\'e}finies par des équations diff{'e}rentielles. 
 1901, Acta Mathematica, Springer Netherlands, 24 (1): 1–88., doi:10.1007/BF02403068.
 

\bibitem{BenGou04-1}
Bennett, D.~L., and Gourley, S.~A. Global stability in a model of the glucose--insulin interaction with time delay. Euro. Jnl of Appl. Math. 2004, 15, 203--221.

\bibitem{BenGou04-2} Bennett, D.~L., and Gourley, S.~A. Periodic oscillations in a model of the glucose--insulin interaction with delay and periodic forcing. Dynamical Systems, 2004, 19(2), 109-125.

\bibitem{BenGou04-3} Bennett, D.~L., and Gourley, S.~A. Asymptotic properties of a delay differential equation model for the interaction of glucose with plasma and interstitial insulin. Applied Mathematics and Computation, 2004, 151, 189-207. 




\bibitem{BerTre2004} Berrut, J.-P., and Trefethen, L.~N. Barycentric Lagrange Interpolation. SIAM Review, 2004, 46, 501.

\bibitem{BraHasIvaTro20} Braverman, E., Hasik, K., Ivanov, A., and Trofimchuk, S. A cyclic system with delay and its characteristic equation. Discrete and Continuous Dynamical Systems, Ser. S. 2020, 13(1), 1-29.

\bibitem{Coddington}  Coddington, E.A., Levinson, N. The Poincar \'{e}–Bendixson Theory of Two-Dimensional
Autonomous Systems. Theory of Ordinary Differential Equations. 1955, New York: McGraw-Hill. 389–403, 
ISBN 978-0-89874-755-3.

\bibitem{deMStr93} de Melo, W., and van Strien, S. One-dimensional dynamics. Ergebnisse der Mathematik und ihrer Grenzgebiete 3 [Results in Mathematics and Related Areas 3]. 1993, vol. 25.
Springer-Verlag, Berlin, 605 pp.


\bibitem{DieSvGSVLWal95} Diekmann, O., van Gils, S., Verdyn Lunel, S., and Walther, H.-O. Delay Equations: Complex, Functional,  and Nonlinear Analysis. 1995, Springer-Verlag, New York.


\bibitem{engelborgs2001} Engelborghs, K., Lemaire, V., B\'elair, J., Roose, D. Numerical bifurcation analysis of delay differential equations arising from physiological modeling. J. Math. Biol. 2001, 42, 361.

\bibitem{HadTom77} Hadeler, K.~P., and Tomiuk, J. Periodic solutions of difference differential equations.
Arch. Rat. Mech. Anal. 1977, 65, 87-95.

\bibitem{HalIva93}
J.K. Hale and A.F. Ivanov,  On a high order differential delay
equation. {\it J. Math. Anal. Appl.} {\bf173} (1993), 505--514.

\bibitem{HalSVL93} Hale, J.~K., and Verduyn Lunel, S.~M. Introduction to Functional Differential Equations. 1993, vol. 99. Springer Applied Mathematical Sciences.

\bibitem{Hansen23} Hansen, K. Oscillations in the blood sugar in fasting normal persons. Acta Med. Scand. Suppl. 1923,  4, 27-58.

\bibitem{HirSma74}
Hirsch, M. W. and Smale, S. Differential Equations, Dynamical Systems, and Linear Algebra. 
Ser.: "Pure and Applied Mathematics," vol. 60, 359 pp. Academic Press, 1974.

\bibitem{HuaBriAng17} Huard, B., Bridgewater, A., and Angelova, M. Mathematical investigation of diabetically impared ultradian oscillations in the glucose-insulin regulation. J. Theor. Biology, 2017, 418, 66-76.

\bibitem{HuaEasAng15}
Huard, B., Easton, J.~F., and Angelova, M. Investigation of stability in a two-delay model of the ultradian oscillations in glucose-insulin regulation. Commun. Nonlinear Sci. Numer. Simulat. 2015, 26, 211-222.


\bibitem{IvaDza20}
Ivanov, A., and Dzalilov, Z.
Oscillations and periodic solutions in a two-dimensional differential delay model.
Proceedings of the international conference AMMCS-2019, Springer-Verlag, 2020, 11 pp. (to appear)

\bibitem{IvaBLW04} Ivanov, A.~F., and Lani-Wayda, B. Periodic solutions for three-dimensional non-monotone cyclic systems with time delays. Discrete and Continuous Dynamical Systems. 2004, 11 (2,3), 667-792.

\bibitem{IvaBLW19} 
Ivanov, A.~F., and Lani-Wayda, B. Periodic solutions for an $N$-dimensional cyclic feedback system with delay. 
J. Differential Equations 2020, 268, 5366--5412.

\bibitem{IvaSha91} Ivanov, A.~F., and Sharkovsky, A.~N. Oscillations in singularly perturbed delay equations. Dynamics Reported (New Series), 1991, 1, 165-224.

\bibitem{keenerandsnayd} Keener, J., Sneyd, J. Mathematical Physiology. 1998, Springer, New York.

\bibitem{kissler2014} Kissler, S., Cichowitz, C., Sankaranarayanan, S., Bortz, D.  Determination of personalized diabetes treatment plans using a two-delay model. J. Theor. Biol. 2014, 359, 101-111.


\bibitem{BLW1} Lani-Wayda, B. Persistence of Poincar{\'e} mappings in functional-differential equations (with application to structural stability of complicated behavior).  J. Dynam. Differential Equations. 1995, 7(1), 1-71.

\bibitem{BLW2} Lani-Wayda, B. Hyperbolic sets, shadowing and persistence for noninvertible mappings in Banach spaces. Pitman Research Notes in Mathematics Series, 1995, 334. Longman, Harlow.


\bibitem{LiKua07} Li, J. and Kuang, Y. Analysis of a model of the glucose-insulin regulatory system with two delays. SIAM J. Appl. Math. 2007, 67, 757-776.


\bibitem{LiKuaMas06} Li, J., Kuang. Y., and Mason, C. Modeling the glucose-insulin regulatory system and ultradian insulin secretory oscillations with two time delays. J. Theoret. Biol. 2006, 242, 722-735.

\bibitem{LiWang2012} Li, J., Wang, M., De Gaetano, A., Palumbo, P., Panunzi, S. The range of time delay and the global stability of the equilibrium for an ivgtt model. Math. Biosci. 2012, 235, 128-137.

\bibitem{M-P88}
 Mallet-Paret, J. Morse decompositions for delay differential
equations. J. Differential Equations 1988, 72  270--315.

\bibitem{M-PSel96a} 
Mallet-Paret, J. and Sell, G.~R. Systems of differential delay equations: Floquet multipliers and discrete Lyapunov functions. 
J. Differential Equations, 1996, 125, 385--440.

\bibitem{M-PSel96b} Mallet-Paret, J. and Sell, G.~R. The Poincar{\'e}-Bendixson theorem for monotone cyclic feedback systems  with  delay. J. Differential Equations, 1996, 125, 441--489.

\bibitem{M-PWal94} Mallet-Paret, J., and Walther, H.~O. Rapid oscillations are rare in scalar systems governed 
by monotone negative feedback with a time lag. Preprint, 1994, 35pp.

\bibitem{marchetti2016} Marchetti, L., Reali, F., Dauriz, M., et al.,  A Novel Insulin/Glucose Model after a Mixed-Meal Test in Patients with Type 1 Diabetes on Insulin Pump Therapy. 2016, Scientific Reports, 6, 36029.

\bibitem{Nus79b} Nussbaum, R.~D. Uniqueness and nonuniqueness of periodic solutions of \newline
$x^\prime(t)=g(x(t-1))$. J. Differential Equations, 1979, 34, 25-54.

\bibitem {Poincare} Poincar{\'e}, H. Sur les courbes d{\'e}finies par une {\'e}quation différentielle. 
1892, Oeuvres, 1, Paris

\bibitem{Satin2015} 
Satin, L.~S., Butler, P.~C., Ha, J., Sherman, A.~S. Pulsatile insulin secretion, impaired glucose tolerance and type 2 diabetes. Mol. Aspects Med., 2015, 42, 61-77.

\bibitem{Scheen96} Scheen, A., Sturis, J., Polonsky, K., Van Cauter, E., Alterations in the ultradian oscillations of insulin secretion and plasma glucose in aging. Diabetologia, 1996, 39(5), 564-572.

\bibitem{ShaKolSivFed97} Sharkovsky, A.~N., Kolyada, S.~F., Sivak, A.~G., and Fedorenko, V.~V. Dynamics of One-dimensional Maps. 1997, 407, 261pp. Kluwer Academic Publishers, Ser.: Mathematics and Its Application.

\bibitem{ShaMaiRom93} Sharkovsky, A.~N., Maistrenko, Yu.~L., and Romanenko, E.~Yu., Difference Equations and Their Perturbations.  1993, Kluwer Academic Publishers, Ser.: Mathematics and Its Application, 250, 358 pp.

\bibitem{Wal81a} Walther, H.-O. Density of slowly oscillating solutions of $\dot x(t)=-f(x(t-1))$. J. Math. Anal. Appl. 1981, 79(1), 127-140.


\end{thebibliography}
\end{document}